\documentclass[12pt]{article}

\usepackage{mathrsfs}
\usepackage{amsmath}
\usepackage{amssymb}
\usepackage{amsbsy}
\usepackage{cite}
\usepackage{amsfonts}
\usepackage{cite}
\usepackage{float}

\usepackage{tikz}
\usepackage{pgflibraryarrows}
\usepackage{pgflibrarysnakes}

\input amssym.def
\input amssym.tex

\newtheorem{lem}{Lemma}[section]%
\newtheorem{thm}[lem]{Theorem}%
\newtheorem{defi}[lem]{Definition}%
\newtheorem{cor}[lem]{Corollary}%
\newtheorem{exam}[lem]{Example}%
\newtheorem{prop}[lem]{Proposition}%
\newtheorem{rem}[lem]{Remark}%

\newcommand{\ZZ}{\mathbb{Z}}

\newcommand{\PSL}{\hbox{\rm PSL}}

\newcommand{\Aut}{\hbox{\rm Aut\,}}
\newcommand{\Inn}{\hbox{\rm Inn}}

\newcommand{\Cay}{\hbox{\rm Cay }}
\newcommand{\CH}{\hbox{\rm CH}}
\newcommand{\CD}{\hbox{\rm CD}}

\def\Rot{\hbox{\rm Rot}}

\newcommand{\Core}{\hbox{\rm Core}}

\newcommand{\prf}{\f {\bf Proof}\hskip10pt}

\newcommand{\PGL}{\hbox{\rm PGL}}

\newcommand{\Sym}{\hbox{\rm Sym}}

\newcommand{\CM}{\hbox{\rm CM}}

\newcommand{\MM}{\mathcal{M}}
\newcommand{\HH}{\mathcal{H}}

\newcommand{\x}{\mathbf{x}}

\def\a{\alpha}  \def\g{\gamma}  \def\e{\varepsilon}
 \def\s{\sigma} \def\t{\tau}

    \def\G{\varGamma}
\def\ol1{\overline 1}
\def\di{\bigm|}
\def\lg{\langle}
\def\rg{\rangle}

\def\f{\noindent}
\def\qed{\hfill $\Box$\vskip 0.3cm}

\newcommand{\Mon}{\hbox{\rm Mon}}

\def\di{\bigm|}
\begin{document}

\title{\vspace{-1.9cm}\bf Cayley hyper-digraphs and Cayley hypermaps
\footnote{Corresponding to
pktide@163.com (Kai Yuan) and wang$_{-}$yan@pku.org.cn (Yan Wang).}}

\author{
Kai Yuan, Yan Wang \\
{ School of Mathematics and Information Science} \\
{ Yantai University, Yantai, China.}
}

\date{}
\maketitle
\vskip 5mm

\begin{abstract} A Cayley hyper-digraph is a directed hypergraph that its automorphism group contains a subgroup acting regularly on vertices and a Cayley hypermap is a hypermap whose automorphism group
 contains a subgroup which induces  regular action on  the hypervertex set. In this paper, we study Cayley hyper-digraphs and construct Cayley hypermaps which have high level of symmetry. Our main goal is to present the general theory so as to
make it clear to study Cayley hypermaps.

\vskip 2mm

\medskip

\noindent{\bf Keywords:} Cayley hyper-digraph; Cayley hypergraph; Simple hypergraph; Cayley hypermap

\noindent{\bf 2000 Mathematics subject classification:}  05C30, 05C65
\end{abstract}

\medskip

\section{Introduction}
Let $\G = (V,E)$ be a  hypergraph with  hypervertex set
$V$, and with  hyperedge set $E$, where $E\subseteq V^P$, the power set of $V$ and $\emptyset\not\in E$.
Every hypergraph   $\G$ is associated with a graph, so called {\it Levi graph}  ${\cal G}({\G})$,
 with vertex set $V\sqcup E$, in which $v \in V$ and $\mathbf{e} \in E$
are adjacent if and only if $v$ and $\mathbf{e}$ are incident in $\G$.
Let $\G = (V,E)$ be a  hypergraph.
Let $v_i, v_j \in V$ and $\mathbf{e} \in E$. If $v_i, v_j\in \mathbf{e}$, then $v_i$ (also $v_j$) is called {\it incident} with $\mathbf{e}$, $v_i$ and $v_j$ are called {\it adjacent},
 and $(v_i,\mathbf{e})$  is called an {\it arc} of $\G$. The number of hyperedges incident with $v_i$ is called the {\it degree} of $v_i$.  If each hyperedge contains $r$ vertices for some positive integer $r$, then this hypergraph is called a {\it $r$-uniform} hypergraph.
  A {\it simple hypergraph} is a  hypergraph  $(V, E)$ such that for any hyperedge $\mathbf{e}$, $\mathbf{e}$ is not a subset of some other hyperedges.
 A simple hypergraph is  vertex transitive if for any pair of vertices $x, y \in V$ there is a hypergraph automorphism $\s$
such that $x^{\s} = y$.

\vskip 3mm
In this paper, we shall introduce two new concepts. By $V^P$ we denote the power set of a set $V$.  Firstly,  we generalize the concept of hypergraphs and introduce  {\it hyper-digraphs}.
\begin{defi} Let $V$ be a finite set and $D$ be a non-empty subset of $V\times V^P$ such that for any $(u,\x)\in D$, we have $u\in \x$.
Then we call   $\HH=(V,D)$  a {\it hyper-digraph}.
\end{defi}

Here we give three remarks: (i)  every  hyper-digraph $\HH$ satisfying that    for any $(u,\x)\in D$,    if $v\in  \x$ then   $(v, \x)\in D$,   is
   a hypergraph; (ii)  every  hyper-digraph $\HH=(V,D)$ corresponds to a hypergraph $\G=(V,E)$, where $E=\{ \x \di (u, \x)\in D, \, {\rm for\, some}\, u\in V\}$, where $\G$ is called the {\it underlying
hypergraph} of $\HH$; and (iii)   if for every pair $(u, \x)\in D$, $\x$ has size 2, then $\HH$ is an   ordinary  digraph. This is the reason we use the terminology of ``hyper-digraphs".
\vskip 3mm

Secondly, we generalize the concept of Cayley digraphs and introduce   {\it Cayley hyper-digraphs}.

\begin{defi}
Let $G$ be a group and  $X$  a set of subsets $\mathbf{x}_1, \mathbf{x}_2, \ldots, \mathbf{x}_d$ of $G$ where  every  $\mathbf{x}_i$ contains the identity element $1_G$ of $G$.
Let
$$V=G \quad {\rm and}\quad D=\{(g,\mathbf{x}g) \di g \in G , \mathbf{x} \in X\}.$$
Then $\HH=(V,D)$  is called a {\it Cayley hyper-digraph}, denoted by  $\HH=\CD(G, X)$,
 where the set  $X$ is called the  {\it Cayley hyperset} of $\cal H$ and the underlying hypergraph $\G$ of $\HH$ is called a  {\it Cayley hypergraph}.  \end{defi}

If $|\mathbf{x}|=r$ for each $\mathbf{x} \in X$, then $\HH$ is  said to be  {\it $r$-uniform}.
Obviously, the 2-uniform Cayley hyper-digraph $\CD(G, X)$ can be viewed as  a Cayley digraph $\Cay(G , S)$, where $S = \bigcup\limits_{i=1}^{d}\mathbf{x}_i\backslash \{1_G\}$, by
identifying every arc $(g,\{g, x_ig\})$ in $D$  with an arc $(g, x_ig)$ in  Cayley digraph, where  $\x_i=\{ 1, x_i\}$.

Cayley hypergraphs are vertex transitive hypergraphs and vertex transitive hypergraphs have been studied in \cite{BC,CL,G,G2,G3,G4,MS,PS,S}.
A topological  problem about hypergraphs is to represent hypergraphs on the surfaces which are called hypermaps. One significant problem in the theory of hypermaps is to construct hypermaps with high level of symmetry and characterize regular hypermaps.
The theory of regular hypermaps is well-developed, and thoroughly
explained in \cite{Co}. An important and convenient way to visualize hypermaps was introduced by Walsh in
\cite{Wa}. The Walsh representation of a hypermap is the embedding of a bipartite graph. Also connections with Galois theory and Grothendieck's programme are described
in \cite{JS2}.

 There have been   some papers dealing with  regular hypermaps with given automorphism groups.
  The classification of regular hypermaps with  automorphism groups
  $\PSL(2, q)$ or $\PGL(2,q)$ can be extracted from \cite{Sa} by Sah.
  Conder, Poto\v{c}nik and \v{S}ir\'{a}\v{n} extended Sah's investigation to reflexible hypermaps, on both nonorientable and orientable surfaces, and provided explicit generating sets for projective linear groups, see \cite{CPS}.
 In \cite{WNH}, Wang, Nedela and Hu investigated the classification of uniquely regular hypermaps with a nilpotent automorphism group. The problem is reduced to the classification of finite maximally automorphic $p$-group that the order of its automorphism group attains Hall's upper bound.

A  primer hypermap means a hypermap whose automorphism group acts faithfully on the hyperfaces.
 The studies of these  hypermaps was initiated by Breda d'Azevedo  and Fernandes in   \cite{BF}, in there  the authors classified the primer hypermaps with a prime number $p$ of hyperfaces.
 Then they determined all regular oriented hypermaps with  $p$ hyperfaces, see \cite{BF1}. In \cite{DH}, Du and Hu  classified primer hypermaps  with a
product of two prime numbers of hyperfaces. Du and Yuan characterized primer hypermaps with nilpotent automorphism groups and prime hypervertex valency, see \cite{DY}.
Recently, Yuan and Wang give a classification of minimal Frobenius hypermaps, see \cite{YW}.

There are other results dealing with hypermaps. For instance,
in \cite{C}, Conder  determined  all regular
hypermaps of genus 2 to 101, and all non-orientable regular hypermaps of genus 3 to 202;
and in \cite{JM}, Jing and Mohar determined the genus and the non-orientable genus of complete
3-uniform hypergraphs $K_n^3$ when $n$ is even. Here  we are not able to list all these results completely.
 For more background one may see \cite{JW}.

 \vskip 3mm

There have been lots of research papers on  Cayley maps, that is embeddings of Cayley graphs into oriented surfaces. Among of them,  in \cite{RSJTW}, Richter et al.  presented a theory of Cayley maps.
This motivates us to study  Cayley hyper-digraphs and  hypermaps with high level of symmetry, see below.

Some properties of Cayley hypergraphs will be  given in Section 2;  Cayley hypermaps will be investigated in Section 3; and   some questions for further
research will be posed in Section 4.

\section{Cayley hyper-digraphs}
In this section,  basic properties of  Cayley hyper-digraphs and Cayley hypergraphs  will be discussed.

\subsection{Basic Properties}
The first  theorem  below mainly shows that every hyper-digraph with a subgroup $G\le \Aut(\HH)$ acting regularly on hypervertices must be isomorphic to a Cayley hyper-digraph. Therefore, we may
study such hyper-digraphs via   Cayley hyper-digraphs.
Let $\HH=(V, D)$  be a hyper-digraph and  $\G=(V,E)$ the corresponding hypergraph, where $E=\{ \x \di (u, \x)\in D, \, {\rm for\, some}\, u\in V\}$.  From now on,
Set  $N_1(v)=\{ e\di E\di v\le e\}$, meaning the neighborhood of $v$ in the incident graph of hypervertices and hyperedge of $\G$.

\begin{thm} \label{thm-di}
Let $\HH=(V, D)$  be a hyper-digraph and  $\G=(V,E)$ the corresponding hypergraph, where $E=\{ \x \di (u, \x)\in D, \, {\rm for\, some}\, u\in V\}$.

Suppose that   a subgroup $G\le \Aut(\HH)$ acts regularly on hypervertices. Pick up an $v\in V$ and write $N_1(v)=\{ e_1, e_2,\ldots , e_d\}$.
For any $i\in [d]$, set $e_i=\{ v^{g_{ij}}\di  1\le j\le |e_i|\}$ and $\x_i=\{ g_{i1}=1, g_{i2}, \ldots , g_{i|e_i|}\}$. Set $X=\{\x_i \di i\in [d]\}$. Then
\begin{enumerate}
\item [\rm (1)]  for any $\e \in D$, we have $\e=(v,e_i)^g$ for some $i\in [d]$ and $g\in G$; and
\item [\rm (2)]  $\HH$ is isomorphic to the Cayley hyper-digraph $\CD(G, X)$.
\end{enumerate}
\end{thm}
\prf
(1) Suppose that $\e=(u,e)$. Since $G$ acts regularly on hypervertices, there exists a unique element $g\in G$ such that $u=v^g$. Thus $\e^{g^{-1}}=(v,e^{g^{-1}})$. And so $e^{g^{-1}}\in N_1(v)$, that is, $e^{g^{-1}}=e_i$ for some $i\in [d]$. Therefore, $\e=(v,e_i)^g$.

(2)  Set $\overline{\HH}=\CD(G, X)$. Let $\s$ be a bijection from $V$ to $G$ by sending  $v^h$ to $h$,  for any $h\in G$. By (1),  for any $(u,e) \in D$, there
exists $i\in [d]$ and $g\in G$ such that $(u,e)=(v,e_i)^g$. Then $$(u,e)^{\s}=(v,e_i)^{g\s}=(v^g,e_i^g)^\s=(g,\x_ig)\in D(\overline{\HH}).$$
For any $j \in [d]$ and
$h \in G$, we have $(h,\x_jh)^{\s^{-1}} = (v^h,e_j^h)=(v,e_j)^h\in D$.
So $\s$ induces a bijection from $D$ to $D(\overline{\HH})$.
Therefore, $\HH$ is isomorphic to $\overline{\HH}$.
\qed

The following facts are basic for Cayley hyper-digraphs.

\begin{prop}\label{basic}
Let $\HH = \CD(G,X)$ be a Cayley hyper-digraph of $G$ with Cayley hyperset $X$. Then
\begin{enumerate}
\item [\rm (1)]
 $\Aut(\HH)$ contains the right regular representation $G_R$ of $G$ so that $\HH$ is vertex-transitive.

 \item [\rm (2)] The underlying hypergraph $\G$ of $\HH$ is  connected  if and only if $G = \lg \bigcup\limits_{\mathbf{x}\in X}\mathbf{x}\rg$.

\item [\rm (3)] $\HH$ is a hypergraph if and only if $X=\{\mathbf{x}g^{-1} \di \mathbf{x}\in X, g\in \mathbf{x}\}$. In this case, $\bigcup\limits_{\mathbf{x}\in X}\mathbf{x}$ is inverse closed.

\end{enumerate}
\end{prop}
\prf
(1) Recall that $G_R=\lg h_R: x\mapsto xh \di x,h\in G\rg\leq \Sym(G)$.  For any arc  $(g,\mathbf{x}g)$ of $\HH$, where $g\in G$ and $\mathbf{x}\in X$, we have
$$(g,\mathbf{x}g)^{h_R}=(g^{h_R},(\mathbf{x}g)^{h_R})=(gh,\mathbf{x}gh)\in D(\HH),$$
and so $h_R \in \Aut(\HH)$. Therefore $\HH$ is vertex-transitive.
\vskip 3mm

(2) Since $\HH$ is vertex transitive, $\G$ is a connected hypergraph if and only if for any vertex $g\in G$ there is a path leading  $1_G$ to $g$, that is,  there are some elements $s_i\in \mathbf{x_i}$ and $\mathbf{x_i}\in X$ for $i\in [n]$ such that
$$1_G,\mathbf{x_1}, s_1,\mathbf{x_2}s_1,s_2s_1,\ldots,\mathbf{x_n}s_{n-1}\cdots s_2s_1, s_ns_{n-1}\cdots s_2s_1=g.$$
This is true if and only if $G = \lg \bigcup\limits_{\mathbf{x}\in X}\mathbf{x}\rg$.

\vskip 0.2cm
(3) Suppose that $\HH$ is a hypergraph. Clearly,  $X\subseteq \{\mathbf{x}g^{-1} \di \mathbf{x}\in X, g\in \mathbf{x}\}$.
 For any $\mathbf{x} \in X$, since $\HH$ is a hypergraph, we have $(g,\mathbf{x})\in D(\HH)$ for each $g\in \mathbf{x}$. So there exists a $\mathbf{y}\in X$ such that $(g,\mathbf{y}g)=(g,\mathbf{x})$. Hence $\mathbf{y}g=\mathbf{x}$, and so $\mathbf{x}g^{-1}=\mathbf{y}\in X$. It follows that $\{\mathbf{x}g^{-1} \di x\in X,g\in \mathbf{x}\}\subseteq X$. Therefore,
 $X=\{\mathbf{x}g^{-1} \di \mathbf{x}\in X, g\in \mathbf{x}\}$.

Suppose that  $X=\{\mathbf{x}g^{-1} \di \mathbf{x}\in X, g\in \mathbf{x}\}$.
Then for every element $\mathbf{x}\in X$ and each $g\in \mathbf{x}$, since $\mathbf{x}g^{-1}\in X$, we have  $(g,\mathbf{x})=(g, (\mathbf{x}g^{-1}) g)$ is an arc. Therefore, $\HH$ is a hypergraph. \qed

\begin{defi}
Let $X$ be a Cayley hyperset of a group $G$. Then  $[X]=\{\mathbf{x}g^{-1} \di \mathbf{x}\in X, g\in \mathbf{x}\}$ is called the {\it Cayley closure} of $X$. If
$X=[X]$, then $X$ is called {\it Cayley closed}. In particular,  $X$ is called {\it single Cayley closed}, provided $X=\{\mathbf{x}g^{-1} \di g\in \mathbf{x}\}$.
\end{defi}

The following two propositions are trivial.
\begin{prop}
Let $X$ be a Cayley hyperset of group $G$. Suppose each element of $X$ is a subgroup of $G$. Then $X$ is Cayley closed  and the Cayley hyper-digraph $\HH = \CD(G,X)$ is a hypergraph. Moreover, $|E(\HH)|=\sum\limits_{\mathbf{x}\in X}^{}|G:\mathbf{x}|$.
\end{prop}

\begin{prop}
Let $G$ be a group and let $X$ be Cayley closed. Then the Cayley hyper-digraph $\HH = \CD(G,X)$ is a hypergraph.
\end{prop}

\begin{defi}
Let $\mathbf{x}$ and $\mathbf{y}$ be two subsets of a group $G$ containing $1_G$. Then $\mathbf{x}$ and $\mathbf{y}$ are called {\it Cayley equivalent} if there exists a $g\in \mathbf{x}$ such that $\mathbf{x}g^{-1}=\mathbf{y}$. Let $Y$ be a Cayley hyperset of $G$. Then $Y$ is called {\it non-Cayley equivalent}  if
any two different elements of $Y$ are not Cayley equivalent, where   $Y=\{\mathbf{y}\}$   is viewed to be non-Cayley equivalent.
\end{defi}

\begin{defi}
Let $X$ be a Cayley hyperset. Then a non-Cayley equivalent subset  $Y$ of $X$ is called  {\it maximal}, if
 for any $\mathbf{x}\in X$ there exists a $\mathbf{y}\in Y$ such that $\mathbf{x}$ and $\mathbf{y}$ are {\it Cayley equivalent}.
\end{defi}

\begin{defi}\label{uc}
Let $G$ be a group and let $Y$ be a non-Cayley equivalent hyperset. By $\G=\CH(G,Y)$ we denote  the hypergraph with vertex set $V(\G)=G$ and edge set $E(\G)=\{\mathbf{x}g \di \mathbf{x}\in Y, g\in G\}$.
\end{defi}

 Note that $G$ induces a natural action on $E(\G)$ by right multiplication.

\begin{rem}{\rm
In the paper \cite{LK}, Lee and Kwon introduced Cayley hypergraphs which are similar
to definition \ref{uc}.
}
\end{rem}

Obviously, we have the following propositions.
\begin{prop}
Let $G$ be a group and let $Y$ be a non-Cayley equivalent hyperset. Let $\G=\CH(G,Y)$. Then $\G$ is the underlying hypergraph of the Cayley hyper-digraph $\HH=\CD(G,[Y])$.
\end{prop}

\begin{prop}
Let $G$ be a group and let a hyperset $X$ be Cayley closed. Suppose that $Y$ is a maximal non-Cayley equivalent subset of $X$. Let $\G=\CH(G,Y)$. Then $\G$ is the underlying hypergraph of the Cayley hyper-digraph $\HH=\CD(G,X)$.
\end{prop}

\begin{prop}\label{estab}
Let $G$ be a group and let $Y$ be a non-Cayley equivalent hyperset. Let $\G=\CH(G,Y)$. Suppose that $G$ acts on $E(\G)$ by right multiplication.
Take
 $\mathbf{y} \in [Y]$. Then the stabilizer $G_{\mathbf{y}}$ is a subset of $\mathbf{y}$ and $\mathbf{y}$ is the union of some left cosets of $G_{\mathbf{y}}$.
\end{prop}

\begin{exam}
Let $\HH = \CD(G , X)$ be a Cayley hyper-digraph, where $G = \ZZ_n, n\geq 7, X = \{\{0,1,3\},\{0,-1,2\},\{0,-2,-3\}\}$. Obviously, $X$ is single Cayley closed and so $\HH$ is a hypergraph. Now $V = \ZZ_n$,
$$D=\{(g, \mathbf{x}+g\di g\in \ZZ_n,  \mathbf{x}\in X\},$$
$$E(\HH) = \{g+\mathbf{x}\di  g\in \ZZ_n,  \mathbf{x}\in X\}=\{\{0, 1, 3\}, \{1,2,4\},\{2,3,5\}, \ldots, \{n-1, 0, 2\}\}.$$
Clearly, if $n=7$, then this hypergraph is the Fano plane, which is  the hypergraph $\CH(G,Y)$ with $Y=\{\{0,1,3\}\}$. If $n=8$, then this hypergraph can be drawn on Projective Plane as in Figure~\ref{fig1}.

\begin{figure}[ht]
 \begin{center}
  \begin{tikzpicture}[
dashed line/.style={dashed, thin},
dot/.style = {
      draw,
      fill = white,
      circle,
      inner sep = 0pt,
      minimum size = 5pt
    }
]
\draw (0,0) -- (2,0) -- (3,0)-- (4,0) -- (6,0) (3,0)--(2,1)--(2,-1)--(3,0) (3,0)--(4,1)--(4,-1)--(3,0) (2,1)--(0,0)--(2,-1) (4,1)--(6,0)--(4,-1);

\draw (6,0) arc (0:360:3 and 2);

\fill (0,0) circle (3pt); \fill(2,0) circle (3pt); \fill(3,0) circle (3pt);\fill(4,0) circle (3pt);\fill(6,0) circle (3pt); \fill(2,1) circle (3pt);\fill(2,-1) circle (3pt);\fill(4,1) circle (3pt);\fill(4,-1) circle (3pt);

\path (0,0) node[label = {left: 6}] {};
\path (6,0) node[label = {right: 6}] {};
\path (2,-1) node[label = {below: 7}] {};
\path (2,1) node[label = {above: 4}] {};
\path (4,-1) node[label = {below: 1}] {};
\path (4,1) node[label = {above: 0}] {};
\path (3,0) node[label = {below: 2}] {};
\path (2.1,0.3) node[label = {left: 5}] {};
\path (3.9,0.3) node[label = {right: 3}] {};

\end{tikzpicture}
\end{center}

 \caption{$\CH(\ZZ_8,Y)$ drawn on Projective Plane with $Y=\{\{0,1,3\}\}$.}
 \label{fig1}
\end{figure}

\end{exam}

\subsection{Automorphisms of Cayley hyper-digraphs}

\vskip 0.5cm
For a Cayley hyper-digraph $\HH = \CD(G , X)$, by Proposition~\ref{basic}, for each $h\in G$, the right translation $h_R: g\mapsto gh$ is an automorphism of the Cayley
hyper-digraph $\HH = \CD(G , X)$. The group $G_R = \{h_R \di h \in G \}$ (isomorphic to $G$) is a subgroup of $\Aut (\HH)$ acting
regularly on $V(\HH)$. We will find out the normalizer of $G_R$ in $\Aut (\HH)$.

Some automorphisms of $G$ may induce automorphisms of the Cayley hyper-digraph $\HH$.  Such automorphisms lie in the following
subgroups:
$$\Aut(G, X) = \{\s\in \Aut(G) \di \mathbf{x}^{\s} \in X, \hbox{~for~any~}\mathbf{x}\in X\},$$
$$\Inn(G, X) = \{\s\in \Inn(G) \di \mathbf{x}^{\s} \in X, \hbox{~for~any~}\mathbf{x}\in X\}.$$

\begin{lem}\label{lastlemma}
For a Cayley hyper-digraph $\HH = \CD(G , X)$, we have
$$\Aut(G) \cap \Aut(\HH) = \Aut(G, X),$$
$$\Inn(G) \cap \Aut(\HH) = \Inn(G, X).$$
\end{lem}

\prf
Let $\s\in \Aut(G) \cap \Aut(\HH) $.  Then $\s$ fixes $1_G$ and  then fixes   all the edges containing the vertex $1_G$, which implies $\s$ fixes $X$  setwise.
 Therefore, $\Aut(G) \cap \Aut(\HH) \le \Aut(G, X)$.

Conversely, suppose  $\s \in \Aut(G, X)$.  Then $X^{\s} =X$.
For any $(g,\mathbf{e})\in D(\HH)$, there exists $\mathbf{x}\in X$ such that
$(g,\mathbf{e})=(1_G,\mathbf{x})^{g_R}$. Then

$$\left.\begin{array}{ccc}
(g,\mathbf{e})\in D(\HH)  &\Longleftrightarrow &(1_G,\mathbf{x})\in D(\HH)~~~~~ \\ &\Longleftrightarrow &\mathbf{x}\in X ~~~~~~~~~~~~~~~
\\ &\Longleftrightarrow &\mathbf{x}^\s\in X~~~~~~~~~~~~~
\\ &\Longleftrightarrow &(1_G,\mathbf{x}^\s)\in D(\HH)~~~
\\ &\Longleftrightarrow &(1_G,\mathbf{x}^\s)^{g^\s_R}\in D(\HH)~
\\ &\Longleftrightarrow &(g^\s,(\mathbf{x}g)^\s)\in D(\HH)
\\ &\Longleftrightarrow &(g,\mathbf{e})^\s\in D(\HH).~~~~
\end{array}\right.$$
\f  Thus $\s \in \Aut(\HH)$, and so $\Aut(G, X) \le \Aut(G) \cap \Aut(\HH)$. Finally,
$$\Inn(G) \cap \Aut(\HH) = \Inn(G) \cap (\Aut(G) \cap \Aut(\HH)) = \Inn(G) \cap \Aut(G, X) = \Inn(G, X).$$
\qed

\vspace{0.2cm}
The following Theorem \ref{thm-2} determines the normalizer of $G_R$ in $\Aut (\HH)$.

\begin{thm}\label{thm-2}
For a Cayley hyper-digraph $\HH = \CD(G , X)$, the normalizer of the regular subgroup $G_R$ in $\Aut (\HH)$
is $G_R\rtimes \Aut(G, X)$.
\end{thm}
\prf
Since the normalizer of $G_R$ in $\Sym(G)$ is $G_R\rtimes \Aut(G)$, we get  from  Lemma~\ref{lastlemma} that
$$\begin{array}{lll}N_{\Aut (\HH)}(G_R)&=&\Aut(\HH) \cap (G_R\rtimes \Aut(G))\\
&=&G_R(\Aut(\HH) \cap\Aut(G))=G_R\rtimes \Aut(G, X).\end{array}$$
\qed

\begin{cor}\label{normalizer}
Let $G$ be a group and let $Y$ be a non-Cayley equivalent hyperset. For a Cayley  hypergraph $\G = \CH(G , Y)$, the normalizer of $G_R$ in $\Aut (\G)$
is $G_R\rtimes \Aut(G, [Y])$.
\end{cor}

\begin{exam}
Let $\HH$ be the Cayley hyper-digraph $\CD(G , X)$, where $G = \ZZ_7$ and $X = \{\{0,1,3\},\{0,4,5\},\{0,2,6\}\}$.
Note that $A=\Aut(G)\cong \ZZ_6$ and $|X|=3$. If $X^{\s}=X$ for some $\s\in A$, then $|\s|=1$ or $|\s|=3$. Clearly, there exists an automorphism $\s$ of order $3$ such that $X^{\s}=X$. Therefore, the normalizer of $G_R$ in $\Aut (\HH)$
is isomorphic to $\ZZ_7\rtimes \ZZ_3$.
\end{exam}

\begin{cor}
Let $G$ be a group and let $S$ be a set of all characteristic subgroups of $G$. Suppose $X\subseteq S$. Then for the Cayley  hypergraph $\HH = \CD(G , X)$, the normalizer of $G_R$ in $\Aut (\HH)$
is $G_R\rtimes \Aut(G)$.
\end{cor}

\begin{cor}
Let $G$ be a group and let $\mathbf{x}$ be a subset of $G$ and $1_G\in \mathbf{x}$. Suppose $X=\{\mathbf{x}^\s \di \s\in \Aut(G)\}$. Then for the Cayley  hyper-digraph $\HH = \CD(G , X)$, the normalizer of $G_R$ in $\Aut (\HH)$
is $G_R\rtimes \Aut(G)$.
\end{cor}

\section{Cayley hypermaps}

In order to investigate Cayley hypermaps, we should introduce ideal cycles.

\subsection{Ideal cycles}
Let $G$ be a group. Suppose that $\mathbf{y}\subseteq G$ with $1_G\in \mathbf{y}$ and $|\mathbf{y}|=m$. Suppose that $\tau \in \Sym_G$ is  a cycle on  $\mathbf{y}$ such that $$\tau=(a_{1}, a_{2}, \ldots, a_{m}).$$
For any $g\in G$,  following  $\tau$   define a cycle $\tau_{g}$ on $\mathbf{y}g$  as
$$\tau_{g}=(a_{1}g, a_{2}g, \cdots, a_{m}g).$$
The cycle $\tau$ is said to be  {\it ideal}, if  $\tau=\tau_{g}$ for any $g\in G$ such that $\mathbf{y}g=\mathbf{y}$. Correspondingly, the subset $\mathbf{y}$ is said to
 be {\it ideal}.
\begin{prop}\label{pro3.1}
Let $G$ be a finite group and $\mathbf{y}=\{1_G, a_2, \ldots, a_m\}$ be an ideal subset of $G$ with an ideal cycle $\t=(1_G, a_2, \ldots, a_m)$. Suppose that $\sigma\in \Aut(G)$ and  $2\leq i\leq m$. Then $(1_G, a_2^\sigma, \ldots, a_m^\sigma)$ and $\t_{a_i^{-1}}$ are ideal cycles on ideal subsets $\mathbf{y}^\sigma$ and $\mathbf{y}a_i^{-1}$, respectively.
\end{prop}
\prf
Set $\tilde \t=(1_G, a_2^\sigma, \ldots, a_m^\sigma)$. Then $\tilde \t$ is a cycle on $\mathbf{y}^\sigma$.
Suppose that  $h\in G$ and $\mathbf{y}^\sigma h=\mathbf{y}^\sigma$. Then
   $\tilde \t_h=(h, a_2^\sigma h, \ldots, a_m^\sigma h)$. The equality $\mathbf{y}^\sigma h=\mathbf{y}^\sigma$ implies
$(\mathbf{y}h^{\sigma^{-1}})^\sigma=\mathbf{y}^\sigma$. Thus $\mathbf{y}h^{\sigma^{-1}}=\mathbf{y}$. Since $\t$ is an ideal cycle, we have
 $$(h^{\sigma^{-1}}, a_2h^{\sigma^{-1}}, \ldots, a_mh^{\sigma^{-1}})=(1_G, a_2, \ldots, a_m),$$
and so $\tilde \t=\tilde \t_h$, meaning $\tilde \t$ is ideal and $\mathbf{y}^\sigma$ is ideal.

Note that $\t_{a_i^{-1}}=(a_i^{-1}, a_2a_i^{-1}, \ldots, a_ma_i^{-1})$ is a cycle on $\mathbf{y}a_i^{-1}$. Suppose that $g\in G$ and
$\mathbf{y}a_i^{-1}g=\mathbf{y}a_i^{-1}$. Then $\mathbf{y}a_i^{-1}ga_i=\mathbf{y}$. Since $\t$ is an ideal cycle, we have
 $$(a_i^{-1}ga_i, a_2a_i^{-1}ga_i, \ldots, a_ma_i^{-1}ga_i)=(1_G, a_2, \ldots, a_m).$$
 It follows that $$(a_i^{-1}g, a_2a_i^{-1}g, \ldots, a_ma_i^{-1}g)=(a_i^{-1}, a_2a_i^{-1}, \ldots, a_ma_i^{-1}).$$
Hence $\t_{a_i^{-1}}$ is an ideal cycle of $\mathbf{y}a_i^{-1}$ and $\mathbf{y}a_i^{-1}$ is ideal.
\qed

\begin{defi}\label{idealhypergraph}
Let $G$ be a group,  $Y=\{\mathbf{y}_1, \mathbf{y}_2, \ldots, \mathbf{y}_r\}$  a non-Cayley equivalent hyperset and $\G = \CH(G , Y)$, a Cayley hypergraph. Then  $\G$ is called an {\it ideal Cayley hypergraph},  if $\mathbf{y}_i$ is ideal for each $i\in [r]$.
\end{defi}

\begin{lem} \label{ideal3}
Let $G$ be a group and $1_G\in \mathbf{y}\subseteq G$. Suppose that $G_{\mathbf{y}} = \{g \di \mathbf{y}g=\mathbf{y} \hbox{~and~} g\in G\}$. Then all  cycles on $\mathbf{y}$ are ideal  if and only if either $G_{\mathbf{y}}=1$ or $|\mathbf{y}|\le 3$.
\end{lem}

\prf Take a cycle $\t $   on  $\mathbf{y}$ such as $${\t }=(a_{1}=1_G, a_{2}, \ldots, a_{s}).$$

 Suppose $G_{\mathbf{y}}=1$. Then $\mathbf{y}g=\mathbf{y}$ holds only if $g=1$ and so $\tau$  is clearly  ideal.   Suppose that $G_{\mathbf{y}}\ne 1$ and $|\mathbf{y}|\le 3$.
 Recall  that, by Proposition \ref{estab},  $\mathbf{y}$  is a union of some left cosets  of  $G_{\mathbf{y}}$ in $G$.  Then  $G_{\mathbf{y}}\cong \ZZ_{|\mathbf{y}|}$, where   $|\mathbf{y}|=2$ or 3.  Check  directly that   $\tau$ is  ideal.

\vskip 3mm
Conversely suppose all cycles on $\mathbf{y}$ are ideal. It suffices that    if $|\mathbf{y}|\ge 4$, then $G_{\mathbf{y}}=1$. Suppose that all cycles on
$\mathbf{y}$ are ideal. We prove it by discuss three cases, separately.

Suppose that  $2 \di |G_{\mathbf{y}}|$.   Then there is an element $g$ of order 2 in $\mathbf{y}$.  By Proposition \ref{estab}, we may write
$${\tau}=(1_G, g, h, hg,\ldots),$$
 for  $h\in \mathbf{y}\setminus \lg g\rg$. Then $\tau_g=(1_Gg=g,g^2=1_G,hg,hg^2=h,\ldots)$. Clearly,
$\tau\ne \tau_g$. And so $\tau$ is not an ideal cycle, a contradiction.

Suppose that  $3 \di |G_{\mathbf{y}}|$. Since $G_{\mathbf{y}}$ is a subgroup,  there is an element $a$ of order 3 in $\mathbf{y}$. By Proposition \ref{estab}, we may assume
$${\tau}=(1_G, a, a^2, h,\ldots),$$
where $h\in \mathbf{y}\setminus \lg a\rg$. Then $\tau_a=(1_Ga=a,a^2,1_G,ha,\ldots)$. Clearly,
$\tau\ne \tau_a$. And so $\tau$ is not an ideal cycle, a contradiction.

Suppose that   there is a prime $p>3$ such that $p \di |G_{\mathbf{y}}|$. then there is an element $b$ of order $p$ in $\mathbf{y}$. By Proposition \ref{estab}, we may assume
$${\tau}=(1_G, b^2, b, b^3,\ldots).$$
Then $\tau_b=(1_Gb=b,b^3,b^2,b^4,\ldots)$. Clearly,
$\tau\ne \tau_b$. And so $\tau$ is not an ideal cycle, a contradiction.

In summary,  $G_{\mathbf{y}}=1$.
\qed

 Lemma \ref{ideal3} will  be used to construct Cayley hypermaps with 3-uniform underlying hypergraphs.

\subsection{Introducing Cayley hypermaps}
Topologically, a {\it hypermap} $\mathcal{M}$ is a $2$-cell embedding  of the Levi graph $\mathcal{G}(\G)$ of a connected hypergraph $\varGamma$ (may have multiple edges) into a compact and connected surface $\mathcal{S}$ without boundary,
where the vertices of $\mathcal{G}(\G)$ in two biparts are respectively called  hypervertices and hyperedges of the hypermap, and the connected regions of $\mathcal{S}\backslash\mathcal{G}(\G)$ are called {\it hyperfaces}. The hypergraph $\varGamma$ is called the {\it underlying
hypergraph} of $\mathcal{M}$.
 A hypermap $\mathcal{M}$ is {\it orientable} if the underlying surface $\mathcal{S}$ is orientable. One may look at \cite{DY} for a detailed definition of a hypermap.

 For a topological hypermap ${\cal M}$, choose a center for each  hyperface and  subdivide the hypermap by adjoining the hyperface centers to its adjacent hypervertices and hyperedges. Then we get a triangular subdivision whose triangles are the {\it flags} of this  hypermap,    which  are represented by little triangle around  hypervertices  as shown in Figure~\ref{flag1} (right).

 \begin{figure}[ht]
 \begin{center}
  \begin{tikzpicture}[
dashed line/.style={dashed, thin},
dot/.style = {
      draw,
      fill = white,
      circle,
      inner sep = 0pt,
      minimum size = 5pt
    }
]
\draw (1,2) -- (2,1) -- (3,2) (1,-2) -- (2,-1) -- (3,-2) (2,1) -- (2,-1);

\draw[dashed](0,0) -- (1,2) (0,0) -- (2,1) (0,0) -- (2,-1) (0,0) -- (1,-2) (4,0) --  (3,2) (4,0) -- (2,1) (4,0) -- (2,-1) (4,0) -- (3,-2) ;

\fill (2,1) circle (3pt); \fill(1,-2) circle (3pt); \fill(3,-2) circle (3pt);

\draw (1,2) node[dot] {} (3,2) node[dot] {} (2,-1) node[dot] {} ;

\path (1,0) node[label = {right: $f^{\gamma_2}$}] {};
\path (3,0) node[label = {left: $f$}] {};
\path (3,1.5) node[label = {below: $f^{\gamma_1}$}] {};
\path (3,-1.5) node[label = {above: $f^{\gamma_0}$}] {};

\node at (2,2) {$\cdots$};

\node at (2,-2) {$\cdots$};
  \begin{scope}[xshift=7cm]

\draw (1,2) -- (2,1) -- (3,2) (1,-2) -- (2,-1) -- (3,-2) (2,1) -- (2,-1);

\fill (2,1) circle (3pt); \fill(1,-2) circle (3pt); \fill(3,-2) circle (3pt);

\draw (1,2) node[dot] {} (3,2) node[dot] {} (2,-1) node[dot] {} ;

\fill (2.26,1.12)--(2.52,1.38)--(2.8,1.12)--cycle;
\fill (1.5,0.44)--(1.9,0.44)--(1.9,0.84)--cycle;
\fill (1.08,-1.82)--(1.08,-1.26)--(1.34,-1.52)--cycle;
\fill (2.94,-1.82)--(2.68,-1.54)--(2.96,-1.26)--cycle;
\fill (2.1,0.44)--(2.1,0.84)--(2.5,0.44)--cycle;
\fill (1.2,1.12)--(1.76,1.12)--(1.5,1.38)--cycle;

\path (1.7,0.42) node[label = {below: $f^{\gamma_2}$}] {};
\path (2.3,0.42) node[label = {below: $f$}] {};
\path (2.62,1.24) node[label = {right: $f^{\gamma_1}$}] {};
\path (2.94,-1.54) node[label = {right: $f^{\gamma_0}$}] {};
\node at (2,2) {$\cdots$};
\node at (2,-2) {$\cdots$};
  \end{scope}
\end{tikzpicture}
\end{center}

 \caption{Flags are drawn in the  left figure and   represented by little triangles around  hypervertices  in the right figure.}
 \label{flag1}
\end{figure}
 Now, we define three involutary permutations $\g_0$, $\g_1$ and $\g_2$  on     the set $F$  of flags of ${\cal M}$ as follows:
 $\gamma_0$  exchanges two flags adjacent to  the same hyperedge and  center but  distinct  hypervertices;
 $\gamma_1$   exchanges two flags adjacent to the same hypervertex and  center but   distinct  hyperedges;  and $\gamma_2$ exchanges two flags  adjacent  to the same hypervertex and hyperedge but  distinct centers,  see Figure~\ref{flag1}(left).
 \begin{figure}[ht]
\centering
\begin{tikzpicture}[
dashed line/.style={dashed, thin},
dot/.style = {
      draw,
      fill = white,
      circle,
      inner sep = 0pt,
      minimum size = 5pt
    }
]

\draw (1,2) -- (2,1) -- (3,2) (1,-2) -- (2,-1) -- (3,-2) (2,1) -- (2,-1);

\fill (2,1) circle (3pt); \fill(1,-2) circle (3pt); \fill(3,-2) circle (3pt);

\draw (1,2) node[dot] {} (3,2) node[dot] {} (2,-1) node[dot] {} ;

\draw (2.26,1.12)--(2.52,1.38)--(2.8,1.12)--cycle;
\draw (1.5,0.44)--(1.9,0.44)--(1.9,0.84)--cycle;
\fill (1.08,-1.82)--(1.08,-1.26)--(1.34,-1.52)--cycle;
\draw (2.94,-1.82)--(2.68,-1.54)--(2.96,-1.26)--cycle;
\fill (2.1,0.44)--(2.1,0.84)--(2.5,0.44)--cycle;
\fill (1.2,1.12)--(1.76,1.12)--(1.5,1.38)--cycle;
\node at (2,2) {$\cdots$};
\node at (2,-2) {$\cdots$};
 \end{tikzpicture}
 \caption{Two orbits of flags are represented by black   little triangles and white little triangles.}
 \label{flag2}
\end{figure}

 Since ${\cal G}(\G)$ is connected, the subgroup  $\lg \g_0, \g_1, \g_2\rg $ of   $S_F$ acts transitively  on $F$. In the orientable case,
the even word subgroup $\langle\gamma_0\gamma_1,\gamma_1\gamma_2\rangle$ of $\lg \g_0, \g_1, \g_2\rg $ acts on $F$  with two orbits, see Figure~\ref{flag2}.
\begin{figure}[ht]
\begin{center}
\begin{tikzpicture}[
->,>=stealth',shorten >=1pt,
dot/.style = {
      draw,
      fill = white,
      circle,
      inner sep = 0pt,
      minimum size = 5pt
    }
]
\draw[-](-1.9,1.08) to (0,0);
\draw[-](0,0)to(1.9,1.08);
\draw[-](0,0) to (0,-2.2);

\fill (0,0) circle (3pt);
\draw (1.9,1.08) node[dot] {} (-1.9,1.08) node[dot] {} (0,-2.2)node[dot] {};

\fill (0,0.2)--(0.4,0.36)--(0.18,0.78)--cycle;
\fill (-0.62,0.12)--(-0.22,-0.08)--(-0.82,-0.24)--cycle;
\fill (0.16,-0.2)--(0.16,-0.66)--(0.56,-0.66)--cycle;

\path (-0.04,0.44) edge  [bend right] node[above]{$\gamma_1\gamma_2$} (-0.46,0.12)
      (-0.5,-0.26)  edge  [bend right] node[below]{$\gamma_1\gamma_2$} (0.1,-0.62)
      (0.44,-0.38) edge  [bend right] node[right]{$\gamma_1\gamma_2$} (0.34,0.28) ;

\begin{scope}[xshift=7cm]

\draw[-](-1.9,1.08) to (0,0);
\draw[-](0,0)to(1.9,1.08);
\draw[-](0,0) to (0,-2.2);

\fill (-1.9,1.08) circle (3pt)  (1.9,1.08)circle (3pt)   (0,-2.2) circle (3pt);
\draw (0,0) node[dot] {} ;

\fill (1.60,0.14)--(1.78,0.80)--(1.36,0.54)--cycle;
\fill (-1.68,1.24)--(-1.29,1.02)--(-1.0,1.42)--cycle;
\fill (-0.16,-2.08)--(-0.16,-1.64)--(-0.62,-1.64)--cycle;

\path (1.84,0.48) edge  [bend left] node[right]{$\gamma_0\gamma_2$} (0.08,-1.98)
     (-1.44,1.40)   edge  [bend left] node[above]{$\gamma_0\gamma_2$} (1.54,1.14)
   (-0.38,-1.96)   edge  [bend left] node[below right]{$\gamma_0\gamma_2$}  (-1.82,0.9);

\end{scope}
\end{tikzpicture}
\end{center}

\vskip 3mm
\begin{center}
\begin{tikzpicture}[
->,>=stealth',shorten >=1pt,
dot/.style = {
      draw,
      fill = white,
      circle,
      inner sep = 0pt,
      minimum size = 5pt
    }
]

\draw (0,0) -- (0,1.86) -- (1.32,2.92) -- (2.62,1.86) -- (2.62,0) -- (1.32,-1) -- cycle;

\fill (0,1.86) circle (3pt)  (2.62,1.86)circle (3pt)   (1.32,-1) circle (3pt);
\draw (0,0) node[dot] {} (1.32,2.92) node[dot] {}  (2.62,0) node[dot] {};

\fill (1.9,2.28)--(2.32,1.94)--(1.58,1.88)--cycle;
\fill (0.1,1.66)--(0.1,1.12)--(0.64,1.12)--cycle;
\fill (1.52,0)--(1.46,-0.78)--(1.84,-0.4)--cycle;

\path (1.86,-0.16) edge  [bend right] node[above left]{$\gamma_0\gamma_1$} (2.1,1.8)
      (1.56,2.1)  edge  [bend right] node[anchor=north west]{$\gamma_0\gamma_1$} (0.42,1.52)
      (0.28,1) edge  [bend right] node[above right]{$\gamma_0\gamma_1$} (1.4,-0.3) ;

\end{tikzpicture}
\end{center}
 \caption{Three figures show  the orbit of $\langle\gamma_1\gamma_2\rangle$,  $\langle\gamma_0\gamma_2\rangle$  and
 $\langle\gamma_0\gamma_1\rangle$ on $F_1$.}
 \label{flag3}
\end{figure}

Each orbit determines an orientation described by the action of the even word subgroup. Fixing an orientation we get an oriented hypermap. Thus an
 orientable hypermap gives rise to a pair of oriented hypermaps $\mathcal{M}$ and $\mathcal{M}^{-}$,
 which are mirror images of each other.
Moreover, given one of two orbits, say  $F_1$, every  orbit of  $\lg \gamma_0\gamma_1 \rg $, $\lg \gamma_1\gamma_2\rg $  and $\lg\g_0\g_2\rg$  on   $F_1$ is respectively  the flags contained in one hyperface,   the flags around one hypervertex and the flags around one hyperedge. See Figure~\ref{flag3}.

 Algebraically, given a finite set  $D$ and   a transitive group  $\langle \rho,\tau\rangle$ on $D$, we may define  an oriented hypermap ${\cal M}=(D; \rho, \tau)$, where
 the orbits of $\lg \rho\rg $, $\lg\tau\rg $ and $\lg \rho\tau\rg $ on $D$ are called hypervertices,    hyperedges and hyperfaces, respectively, with incidence given by non-empty intersection.  Moreover,   $D$ is called the arc set and the group $\Mon(\mathcal{M})=\langle \rho, \tau\rangle$ is  called the monodromy group of the hypermap. In the case $(\rho\tau)^2=1$, ${\cal M}$ is an oriented map.

\vskip 3mm
Let $G$ be a group,  $Y=\{\mathbf{y}_1, \mathbf{y}_2, \ldots, \mathbf{y}_r\}$  a non-Cayley equivalent hyperset and $\G = \CH(G , Y)$
 the  Cayley hypergraph. Assume
 $$\mathbf{y}_i=\{a_{i1}=1_G, a_{i2}, \ldots, a_{is_i}\}, 1\leq i\leq r$$ and the Cayley closure  $$[Y]=\{\mathbf{y}_ia_{ij}^{-1} \ |\  1\leq i\leq r, \  1\leq j\leq s_i\}=
 \{\mathbf{x}_{1}, \mathbf{x}_{2}, \ldots, \mathbf{x}_{d}\}.$$ Then, the arc  set of $\G$ is  $$D(\G)=\{(g, \mathbf{x}_{i}g)\ |\  1\leq i\leq d, g\in G\}=\{(ah, \mathbf{y}_{j}h)\ |\  1\leq j\leq r, a\in \mathbf{y}_{j},  h\in G\}.$$
Note that one partite set of $\mathcal{G}(\G)$ is $G$, the other partite set is $\{\mathbf{x}_{i}g\ |\ 1\leq i\leq d, g\in G\}$ and in $\mathcal{G}(\G)$, $h$ is adjacent with $\mathbf{x}_{i}h, 1\leq i\leq d$, for each $h\in G$.

 Take a cycle $\tau_i$ of hypervertices in (adjacent with)  $\mathbf{y}_i$ for each $ 1\leq i\leq r$, say $${\tau_i}=(a_{i1}=1_G, a_{i2}, \cdots, a_{is_i}),$$ hereafter we write $a_{ij}^{\tau_i}=a_{i j+1}$, where $1\leq j\leq s_i$ and the subscript addition is modular $s_i$.
For any $g\in G$, we define a cycle $\tau_{ig}$ on $\mathbf{y}_ig$ corresponding to $\tau_i$ as
$$\tau_{ig}=(a_{i1}g, a_{i2}g, \cdots, a_{is_i}g).$$
If for any $g\in G_{\mathbf{y}_i}$, that is, $\mathbf{y}_ig=\mathbf{y}_i$, we have $\tau_i=\tau_{ig}$, then $\tau_i$ is an  ideal cycle.

\vskip 3mm
To embed $\mathcal{G}(\G)$ to a orientable surface, we need to introduce two local rotations on  hyperedges and hypervertices, respectively.

\vskip 3mm
(1)  Local rotation on hyperedges:
\vskip 3mm

Suppose  $\t_1$ is  an ideal cycle, for each $1\leq i\leq r$.
Following the permutation of hypervertices in $\mathbf{y}_i$,  the cycle of arcs that share the same hyperedge  $ \mathbf{y}_i$
  is  $$((1_G, \mathbf{y}_i), (a_{i2}, \mathbf{y}_i), \cdots, (a_{is_i}, \mathbf{y}_i)).$$
For each $g\in G$, the cycle of arcs sharing the same hyperedge  $\mathbf{y}_ig$ is  $$((g, \mathbf{y}_ig), (a_{i2}g, \mathbf{y}_ig), \cdots, (a_{is_i}g, \mathbf{y}_ig)).$$ In this way, a permutation $\tau$ of $D(\G )$ can be defined as $$(ag, \mathbf{y}_ig)^\tau=(a^{\tau_i}g, \mathbf{y}_ig)$$
 for each $g\in G, a\in \mathbf{y}_i, 1\leq i\leq r$.  Since $\tau_i$ is an ideal cycle for $1\leq i\leq r$, $\tau$ is well defined. In fact, assume $(ag, \mathbf{y}_ig)=(bh, \mathbf{y}_ih)$ for some $b\in \mathbf{y}_i$, then $(bh, \mathbf{y}_ih)^\tau=(b^{\tau_i}h, \mathbf{y}_ih)$. The equalities $ag=bh$ and $\mathbf{y}_ig=\mathbf{y}_ih$ imply that $gh^{-1}\in \mathbf{y}_i$ and $b=agh^{-1}$. If $a^{\tau_i}=c$, then $b^{\tau_i}=cgh^{-1}$ because $\tau_i$ is ideal. Thus, $a^{\tau_i}g=b^{\tau_i}h$ implying $\tau$ is well defined.
 Clearly, $\tau$ is determined by $\tau_{[r]}=\{\tau_i \di i \in [r]\}$ and the cycles of $\tau$ are in one-to-one correspondence with the hyperedges of $\G$.

\vskip 3mm
(2)  Local rotation on  hypervertices:
\vskip 3mm

Set a cycle of hyperedges in $[Y]$, say $\rho_{[Y]}=(\mathbf{x}_{1},  \mathbf{x}_{2},  \ldots,  \mathbf{x}_{d})$, hereafter we set $$\mathbf{x}_{i}^{\rho_{[Y]}}=\mathbf{x}_{i+1} \hbox{~for~}  i\in [d-1] \hbox{~and~} \mathbf{x}_{d}^{\rho_{[Y]}}=\mathbf{x}_{1}.$$ Then
   define a permutation $\rho$ of arcs in $D(\G)$ as $(g, \mathbf{x}_ig)^\rho=(g, \mathbf{x}_{i}^{\rho_{[Y]}}g)$
    for all $1\leq i\leq d$ and $ g\in G$.
 So, the cycles of $\rho$ are in one-to-one correspondence with the hypervertices of $\G$. Now we are ready to give the following definition of Cayley hypermaps.

\begin{defi}\label{Cayleyhmap}
Let $G$ be a group. Suppose that $Y=\{\mathbf{y}_1, \mathbf{y}_2, \ldots, \mathbf{y}_r\}$ is a non-Cayley equivalent hyperset of $G$, $[Y]=\{\mathbf{x}_{1}, \mathbf{x}_{2}, \ldots, \mathbf{x}_{d}\}$ and $\G = \CH(G , Y)$
is an ideal Cayley hypergraph.
   Let $\rho_{[Y]}$ be a cycle of hyperedges in $[Y]$ and let $\tau_i$ be an ideal cycle of  $\mathbf{y}_i, 1\leq i\leq r$ and $\tau_{[r]}=\{\tau_i \di i \in [r]\}$. Define two permutations
  $\rho$ and $ \tau$  of $D(\G)$ such that for each $(g, \mathbf{x}_ig)=(ah, \mathbf{y}_jh)\in D(\G)$,
  $$(g, \mathbf{x}_ig)^\rho=(g, \mathbf{x}_i^{ \rho_{[Y]}}g),$$ $$(g, \mathbf{x}_ig)^\tau=(ah, \mathbf{y}_jh)^\tau=(a^{\tau_j}h, \mathbf{y}_jh), $$ where $g, h\in G, a\in \mathbf{y}_j, i\in [d]$ and $j\in [r]$.
   Then  $(D(\G); \rho, \tau)$, also denoted by $\CM(G, Y, \rho_{[Y]}, \tau_{[r]})$, is called a Cayley hypermap of $G$ with underlying Cayley hypergraph $\G$.
\end{defi}

\begin{rem}
In the paper \cite{LK}, the authors also defined a Cayley hypermap with
a little differennt way. They only considered a Cayley hypermap under special
condition that for any $i, j \in [r]$ and for any $a \in \mathbf{y}_i, b \in \mathbf{y}_j$, $\t_i(a)a^{-1} = \t_j (b)b^{-1}$
if and only if $i = j$ and $a = b$. With the special condition, our definition of
Cayley hypermap is the same with theirs.
\end{rem}

It is obvious that  $\CM(G, Y, \rho_{[Y]}, \tau_{[r]})$ is an orientable embedding of $\mathcal{G}(\G)$.

\begin{exam}\label{examfano}
Let $G=\mathbb{Z}_7$, $\mathbf{y}=\{0, 1, 3\}$ and $Y=\{\mathbf{y}\}$. Then $$[Y]=\{\{0, 1, 3\}, \{0, 4, 5\}, \{0, 2, 6\}\} $$ and
the Cayley hypergraph $\G = \CH(\mathbb{Z}_7, Y)$ is the underlying hypergraph of the Fano plane.
Here $$D(\G)=\{(i, \{i, i+1, i+3\}), (i, \{i, i+4, i+5\}), (i, \{i, i+2, i+6\})\ |\ 0\leq i\leq 6\},$$ and the addition is modular $7$.

Set $\tau_1=(0, 1, 3)$. By Lemma {\rm \ref{ideal3}}, $\tau_1$ is an ideal cycle. Then, $(i, \{i, i+1, i+3\})^\tau=(i+1,  \{i, i+1, i+3\}), (i, \{i, i+4, i+5\})^\tau=(i+4, \{i, i+4, i+5\}),
(i, \{i, i+2, i+6\})^\tau=(i+2, \{i, i+2, i+6\}), 0\leq i\leq 6.$

\noindent {\bf {Case 1.}}
 Set $\rho_{[Y]}=( \{0, 1, 3\}, \{0, 2, 6\}, \{0, 4, 5\})$. Then,
 $(i, \{i, i+1, i+3\})^\rho=(i,  \{i, i+2, i+6\}), (i, \{i, i+2, i+6\})^\rho=(i, \{i, i+4, i+5\}),
(i, \{i, i+4, i+5\})^\rho=(i, \{i, i+1, i+3\}), 0\leq i\leq 6.$

{\rm
In this case, $\mathcal{G}(\G)$ is embedded on Torus. There are $7$ hypervertices, $7$ hyperedges, $7$ hyperfaces (the  hexagons) and $21$ arcs,
see Figure~\ref{fano32} (left). Here,  dashed sides   labelled by common letters of the outside enneagon are pasted together to get the Torus. }

\begin{figure}[H]
 \begin{center}
  \begin{tikzpicture}[
dashed line/.style={dashed, thin},
dot/.style = {
      draw,
      fill = white,
      circle,
      inner sep = 0pt,
      minimum size = 5pt
    }
]

\pgfmathsetmacro{\n}{6} 
  \pgfmathsetmacro{\radius}{3.5} 

  \foreach \i in {1,...,\n} {
    \coordinate (P\i) at ({360/\n * (\i - 1)}:\radius);
  }

  \foreach \i in {1,...,5} {
    \draw[dashed] (P\i) -- (P\the\numexpr\i+1\relax);
  }
  \fill (P2) circle (3pt);
  \fill (P4) circle (3pt);
  \fill (P6) circle (3pt);
  \draw[dashed](P1)--(P6);

  \path (P2) node[label = {above: {\scriptsize $5$}}] {};
  \path (P4) node[label = {left: {\scriptsize $5$}}] {};
  \path (P6) node[label = {below: {\scriptsize $5$}}] {};

  \draw (P4)--(-2,0);
  \draw ({360/12 * 5}:3.03)--(-2,1);
  \draw (P2)--(1,1.55);
  \draw ({360/12 * 3}:3.03)--(-1,1.55);
  \draw ({360/12 }:3.03)--(2,1);
  \draw ({-360/12 }:3.03)--(2,0);
  \draw ({-360/12*2 }:\radius)--(1,-1.55);
    \draw ({-360/12*3 }:3.03)--(0,-2.55);
     \draw ({-360/12*5 }:3.03)--(-1,-1.48);
  \path ({-360/12*3 }:3.03) node[label = {below: {\scriptsize $a$}}] {};
    \path ({-360/12*5 }:3.03) node[label = {left: {\scriptsize $c$}}] {};
    \path ({360/12*5 }:3.03) node[label = {left: {\scriptsize $b$}}] {};
       \path ({360/12*3 }:3.03) node[label = {above: {\scriptsize $a$}}] {};
   \path ({360/12 }:3.03) node[label = {right: {\scriptsize $c$}}] {};
     \path ({-360/12 }:3.03) node[label = {right: {\scriptsize $b$}}] {};

\begin{scope}[xshift=-2cm]
\draw (1,1.5) node[dot] {} (3,1.5) node[dot] {};
\fill (0,1) circle (3pt); \fill (2,1) circle (3pt); \fill (4,1) circle (3pt); 
\draw (0,0) node[dot] {} (2,0) node[dot] {}  (4,0) node[dot] {};
\fill (1,-0.5) circle (3pt);\fill (3,-0.5) circle (3pt);

\path (2.8,1.58) node[label = {right: {\scriptsize $\{0,4,5\}$}}] {};
\path (0.8,1.58) node[label = {right: {\scriptsize $\{0,1,3\}$}}] {};

\draw(1,-0.5)--(1,-1.44) (1.02,-1.54)--(2,-2)--(2.98,-1.52) (3,-1.44)--(3,-0.5);

 \fill (2,-2) circle (3pt);
\draw (3,-1.5) node[dot] {} (1,-1.5) node[dot] {};

\draw(2,-2)--(2,-2.65);

\draw (0.06,-0.06) -- (1,-0.5) -- (1.96,-0.08) (2,0.09)-- (2,1) -- (1.06,1.44) (0.94,1.44)--(0,1)--(0,0.09);
\draw (2.04,-0.08) -- (3,-0.5) -- (3.96,-0.08) (4,0.09)-- (4,1) -- (3.06,1.44) (2.94,1.44)--(2,1);

\path (1,-0.6) node[label = {above: {\scriptsize 2}}] {};
\path (-0.1,0.9) node[label = {right: {\scriptsize 3}}] {};
\path (1.9,0.9) node[label = {right: {\scriptsize 0}}] {};
\path (4.1,0.9) node[label = {left: {\scriptsize 4}}] {};
\path (3,-0.6) node[label = {above: {\scriptsize 6}}] {};
\path (-0.2,0.2) node[label = {right: {\scriptsize $\{2,3,5\}$}}] {};
\path (1.8,0.2) node[label = {right: {\scriptsize $\{0,2,6\}$}}] {};
\path (3.8,0.2) node[label = {right: {\scriptsize $\{3,4,6\}$}}] {};
\path (0.8,-1.3) node[label = {right: {\scriptsize $\{1,2,4\}$}}] {};
\path (2,-2.1) node[label = {above: {\scriptsize 1}}] {};
\path (2.8,-1.3) node[label = {right: {\scriptsize $\{1,5,6\}$}}] {};

\end{scope}
\begin{scope}[xshift=4.5cm]

\draw (0.1,0) -- (0.8,0);
\draw (6.3,0) -- (5.6,0);


\draw (5.6,0) arc (0:360:2.4);

\fill  (0.8,0) circle (3pt);
\draw (5.6, 0) node[dot] {};

\draw (1.0, 0.95) node[dot] {};
\draw (4.8, 1.8) node[dot] {} (4.8, -1.8) node[dot] {};
\draw (1.0, -0.95) node[dot] {};
\fill (5.4,0.95) circle (3pt) (5.4,-0.95) circle (3pt);
\fill (1.6,1.8) circle (3pt) (1.6,-1.8) circle (3pt);
\fill (3.8,2.33) circle (3pt) (3.8,-2.33) circle (3pt);
\draw (2.6, 2.33) node[dot] {} (2.6, -2.33) node[dot] {};

\draw[dashed] (3.2,3.2)--(1.75,2.83)--(0.59,1.80)--(0.1,0.65)--(0.1,-0.65)--(0.59,-1.80)--(1.75,-2.83)--(3.2,-3.2) (3.2,3.2)--(4.65,2.83)--(5.81,1.80)--(6.3,0.65)--(6.3,-0.65)--(5.81,-1.80)--(4.65,-2.83)--(3.2,-3.2);
\path (0.2,0) node[label = {left: {\scriptsize $f$}}] {} (0.8,0) node[label = {right: {\scriptsize $6$}}] {} ;
\path (6.2,0) node[label = {right: {\scriptsize $e$}}] {} (5.7,0) node[label = {left: {\scriptsize $\{1,2,4\}$}}] {};
\draw (2.58,2.41)--(2.45,3) (2.58,-2.41)--(2.45,-3) (3.82,2.41)--(3.95,3) (3.82,-2.41)--(3.95,-3);

\draw (4.86,1.85)--(5.25,2.3) (4.86,-1.85)--(5.25,-2.3) (1.15,-2.3)--(1.56,-1.85) (1.15,2.3)--(1.56,1.85);

\path (2.58,2.31) node[label = {below: {\scriptsize $\{0,1,3\}$}}] {} (2.45,3) node[label = {above: {\scriptsize $c$}}] {} ;
\path (2.58,-2.31) node[label = {above: {\scriptsize $\{2,3,5\}$}}] {} (2.45,-2.9) node[label = {below: {\scriptsize $g$}}] {} ;

\path (3.82,-2.41) node[label = {above: {\scriptsize $5$}}] {} (3.95,-2.9) node[label = {below: {\scriptsize $d$}}] {} ;

\path (3.82,2.41) node[label = {below: {\scriptsize $0$}}] {} (3.95,3) node[label = {above: {\scriptsize $a$}}] {} ;

\path (1.25,2.5) node[label = {left: {\scriptsize $g$}}] {} (1.46,1.55) node[label = {right: {\scriptsize $3$}}] {} ;
\path (1.25,-2.5) node[label = {left: {\scriptsize $e$}}] {} (1.46,-1.55) node[label = {right: {\scriptsize $2$}}] {} ;

\path (5.05,-2.5) node[label = {right: {\scriptsize $f$}}] {} (4.94,-1.55) node[label = {left: {\scriptsize $\{1,5,6\}$}}] {} ;
\path (5.05,2.5) node[label = {right: {\scriptsize $d$}}] {} (4.94,1.55) node[label = {left: {\scriptsize $\{0,4,5\}$}}] {} ;

\draw (0.92,0.98)--(0.34,1.2) (0.92,-0.98)--(0.34,-1.2) (5.48,0.98)--(6.07,1.2) (5.48,-0.98)--(6.07,-1.2);
\path (0.36,1.3) node[label = {left: {\scriptsize $b$}}] {} (0.84,0.8) node[label = {right: {\scriptsize $\{3,4,6\}$}}] {} ;

\path (0.36,-1.3) node[label = {left: {\scriptsize $a$}}] {} (0.84,-0.7) node[label = {right: {\scriptsize $\{0,2,6\}$}}] {} ;

\path (6.04,-1.3) node[label = {right: {\scriptsize $c$}}] {} (4.76,-0.7) node[label = {right: {\scriptsize $1$}}] {} ;
\path (6.04,1.3) node[label = {right: {\scriptsize $b$}}] {} (4.76,0.7) node[label = {right: {\scriptsize $4$}}] {} ;
\end{scope}
\end{tikzpicture}
\end{center}

 \caption{$\CH(\ZZ_7,Y)$ with $Y=\{\{0,1,3\}\}$ embedded on Torus (the left figure) and the orientable surface of genus $3$ (the right figure).}
 \label{fano32}
\end{figure}

\noindent{\bf {Case 2.}}
Set $\rho_{[Y]}=( \{0, 1, 3\}, \{0, 4, 5\}, \{0, 2, 6\})$.
Then, $(i, \{i, i+1, i+3\})^\rho=(i,  \{i, i+4, i+5\}), (i, \{i, i+4, i+5\})^\rho=(i, \{i, i+2, i+6\}),
(i, \{i, i+2, i+6\})^\rho=(i, \{i, i+1, i+3\}), 0\leq i\leq 6.$

{\rm In this case, $\mathcal{G}(\G)$ is embedded on the orientable surface of genus $3$. There are $7$ hypervertices, $7$ hyperedges, $3$ hyperfaces (the  tetradecagons) and $21$ arcs, see Figure~\ref{fano32} (right). Here,  dashed sides labelled by common letters of the outside tetradecagon  are pasted together to get the orientable surface of genus $3$.}

\end{exam}

\subsection{Automorphism group of Cayley hypermaps}
For  oriented hypermaps $\mathcal{M}=(D;\rho, \tau)$ and $\mathcal{M}'=(D';\rho', \tau')$, an isomorphism  between $ \mathcal{M}$ and $\mathcal{M}'$ is a bijection $\psi:D\rightarrow D'$ satisfying $\rho\psi=\psi \rho'$ and $\tau\psi=\psi \tau'$. This gives the automorphism of $ \mathcal{M}$ when $ \mathcal{M}=\mathcal{M}'$. Use $\Aut(\mathcal{M})$ to denote the  automorphism group of $ \mathcal{M}$. It is straightforward that $|\Aut(\mathcal{M})|\leq|D|$.
An oriented hypermap is called {\it regular} if the action of $\Mon(\mathcal{M})=\langle \rho, \tau\rangle$ on $D$ is regular.
 In this case,  the set $D$ can be replaced by $H:=\Mon(\mathcal{M})$,
 so that $\Mon({\cal M})$ and $\Aut(\mathcal{M})$ can be viewed as  the  left and right regular multiplications of $H$, respectively.
 An oriented regular hypermap ${\cal M}$  can therefore be denoted by a triple
 $\mathcal{M} = (H; a, b )$, where $H=\lg a, b\rg $, $a$ corresponds to $\rho$ above and $b$ to $\tau$.
 Then  the hypervertices (resp. hyperedges and hyperfaces) correspond to right cosets of  $\lg a\rg $, (resp. $\lg b\rg $ and $\lg ab \rg$).
 Given a group $H$, $(H; a_1, b_1)$ is isomorphic to $ (H; a_2, b_2)$ if and only if there exists an automorphism $\s $ of $H$ such that
 $a_1^\s= a_2$ and $b_1^\s =b_2$.

 Let $\Gamma$ be the underlying hypergraph of $\mathcal{M}=(D;\rho, \tau)$ and assume that $\Gamma$ is a simple hypergraph. Take $\sigma\in \Aut(\mathcal{M})$, then $\sigma$ induces a permutation, denoted by $\underline{\sigma}$, on $V(\G)$. To be more clear, because  the cycles of $\rho$
are  in one-to-one correspondence with the hypervertices,  we use $\alpha^{\langle \rho\rangle}=(\alpha, \alpha^\rho, \alpha^{\rho^2}, \ldots, \alpha^{\rho^{\ell-1}}), $ where $\alpha\in D$ and $\ell$ is the length of this cycle, to denote the hypervertex incident to $\alpha^{\rho^{i}}, 0\leq i\leq \ell-1$, and define $(\alpha^{\langle \rho\rangle})^{\underline{\sigma}}=(\alpha^\sigma)^{\langle \rho\rangle}$. Similarly, each hyperedge corresponds to a cycle of $\tau$, denoted by $\alpha^{\langle \tau\rangle}, \alpha\in D$, and $(\alpha^{\langle \tau\rangle})^{\underline{\sigma}}=(\alpha^\sigma)^{\langle \tau\rangle}$. The fact $\sigma\in \Aut(\mathcal{M})$ assures that $\underline{\sigma}$ is well-defined and also a bijection on hypervertices. It is clear that $\underline{\sigma} $ preserves the incidence structure between hypervertices and hyperedges of $\G$. So, $\underline{\sigma}\in \Aut(\G) $.
We claim that under the condition of $\G$ being simple $\Aut(\mathcal{M})$ induces a subgroup of $\Aut(\G)$ that acts  faithfully on $V(\G)$, that is, the kernel of this action is trivial.
\begin{lem}\label{faithful}
Let  $\mathcal{M}=(D; \rho, \tau)$ be a hypermap, and $\G$ be the underlying hypergraph which is simple. Then, $\Aut(\mathcal{M})$ induces a subgroup  $\underline{\Aut(\mathcal{M})}$ of $\Aut(\G)$, which acts faithfully on $V(\G)$.
\end{lem}
\prf
Let $\underline{\sigma_1}, \underline{\sigma_2}$ be induced permutations on $V(\G)$ by $\sigma_1, \sigma_2\in \Aut(\mathcal{M})$, respectively. Since $\sigma_i, i=1,2,$  commutes with both $\rho$ and $\tau$, we have $\underline{\sigma_1}\  \underline{\sigma_2} \in \Aut(\G)$. It follows that $\underline{\Aut(\mathcal{M})}=\{\underline{\sigma}\ |\ \sigma\in \Aut(\mathcal{M})\}$ is a subgroup of $\Aut(\G)$.

Suppose $\underline{\sigma}$ is  induced by $\sigma\in \Aut(\mathcal{M})$ which fixes each hypervertex of $\G$, then $\sigma$ fixes each orbit of $\rho$.
Take any hyperedge $e$ with maximal valency. Because each hypervertex incident to $e$ is fixed by $\underline{\sigma}$, $e^{\underline{\sigma}}$ is a hyperedge that is incident to  each hypervertex of $e$. Considering the choice of $e$ which contains the most number of hypervertices, $e^{\underline{\sigma}}$ and $e$ are parallel hyperedges. This is a contraction to the simplicity assumption of $\G$. Therefore, $e^{\underline{\sigma}}=e$ and consequently $\sigma$ fixes arcs  incident with $e$. By backward induction on the valency of hyperedges,  $\sigma$ can be shown to fix each arc of $\G$. So, $\underline{\Aut(\mathcal{M})}$ acts on $V(\G)$ faithfully.
\qed

\vspace{.3cm}
Note that if $\mathcal{M} = (H; a, b )$ is an oriented regular hypermap, then the hypervertices, hyperedges and hyperfaces are in one-to-one correspondence with the right cosets of $\langle a\rangle, \langle b\rangle$ and $\langle ab\rangle$, respectively.  For a $m$-uniform hypergraph, if every $m$ adjacent hypervertices are
 incident with exactly $t$ distinct hyperedges, we will say that the hypergraph has {\it hyperedge-parallelity} $t$. Specially, the hypergraph is  simple if $t = 1$.
\begin{lem}\label{parallel}
Let $\mathcal{M}=(H; a, b)$ be an oriented  regular hypermap.
Then the hyperedge-parallelity of the underlying hypergraph is $|\bigcap_{x\in \langle b\rangle} \langle b\rangle \langle a\rangle^ x : \langle b\rangle|$, where $\langle a\rangle^ x=x^{-1}\langle a\rangle x$ denotes the conjugation of $\langle a\rangle$ by $x$.
\end{lem}
\prf
Because $\mathcal{M}$ is regular, without any loss of generality, assume $\langle b\rangle$ and $\langle b\rangle h$ are parallel hyperedges for some $h\in H$. The hypervertices in $\langle b\rangle$ are $[\langle a\rangle\langle b\rangle : \langle a\rangle]$ which is the set
of right cosets of $\langle a\rangle$ contained in $\langle a\rangle\langle b\rangle$. So, for all $x\in \langle b\rangle$, $\langle b\rangle h \cap \langle a\rangle x\neq \emptyset$ is equivalent to $h\in \bigcap_{x\in \langle b\rangle} \langle b\rangle x^{-1}\langle a\rangle x$. As a result, there are $|\bigcap_{x\in \langle b\rangle} \langle b\rangle \langle a\rangle^ x : \langle b\rangle|$ hyperedges parallel to $\langle b\rangle$.
\qed

A subgroup $B$ of $A$ is called {\em core free} if $\bigcap_{a\in A} a^{-1}Ba=\{1_A\}$. If $B$ is core free, then $A$ acts faithfully on  the right cosets of $B$ by right multiplication.

 \begin{thm}\label{corefree}
Let $\mathcal{M}=(H; a, b)$ be an oriented  regular hypermap with underlying hypergraph $\G$. If $\G$ is simple, then $\Core_H\langle a\rangle=\{1_H\}$.
   \end{thm}
\prf
 If $\G$ is simple, then according to Lemma~\ref{parallel},
 $|\bigcap_{x\in \langle b\rangle} \langle b\rangle \langle a\rangle^x: \langle b\rangle|=1.$
It is clear that $|\langle b\rangle\Core_H\langle a\rangle  : \langle b\rangle|\leq |\bigcap_{x\in \langle b\rangle} \langle b\rangle \langle a\rangle^x : \langle b\rangle|=1$. So, $\Core_H\langle a\rangle\leq \langle b\rangle$. Furthermore,
$|\langle a\rangle  \langle b\rangle : \langle a\rangle |$ is the number of hypervertices contained in $\langle b\rangle$ which equals to the order of $b$. By $|\langle b\rangle : \langle b\rangle \cap \langle a\rangle|= | \langle b\rangle \langle a\rangle: \langle a\rangle |=|b|$, we have $\langle b\rangle \cap \langle a\rangle=\{1_H\}$, and so $\Core_H\langle a\rangle=\{1_H\}$.
\qed

\vspace{.3cm}
\noindent {\bf Remark.} The result of Theorem~\ref{corefree} is a corollary of Lemma~\ref{faithful} as well. In fact, when $\mathcal{M}=(H; a, b)$ is an oriented regular hypermap, $\Aut(\mathcal{M})\simeq H$. According to Lemma~\ref{faithful}, $H$ acts faithfully on the right cosets of $\langle a\rangle$. Therefore, $\Core_H\langle a\rangle=\{1_H\}$.

Look at a Cayley hypermap $\mathcal{M}=\CM(G, Y, \rho_{[Y]}, \tau_{[r]})$. It is known that $G_R\leq \Aut(\CH(G, Y))$. In fact, $G_R$ is an automorphism group of $\Aut(\mathcal{M})$ as well.

\begin{lem}\label{rightregular}
Let $\mathcal{M}=\CM(G, Y, \rho_{[Y]}, \tau_{[r]})$ be a Cayley hypermap. Then $G_R\leq \Aut(\mathcal{M})$.
\end{lem}
\prf
Take an arc $(g, \mathbf{x}_{i}g)$, where $ g\in G, \mathbf{x}_i \in [Y]$ and $\mathbf{x}_i=\mathbf{y}_jk^{-1}$ for $\mathbf{y}_j\in Y, k\in \mathbf{y}_j  $.
 For each $h\in G$, the right translation $h_R$ induces an action, still denoted by $h_R$, on arcs  as  $ (g, \mathbf{x}_{i}g)^{h_R}=(gh, \mathbf{x}_{i}{gh})$. So, $h_R$ is a permutation of $D(\G)$.

  Because $ (g, \mathbf{x}_{i}g)^{h_R\rho}=(gh, \mathbf{x}_{i}{gh})^\rho=(gh, \mathbf{x}_{i}^{\rho_{[Y]}}gh)$ and
 $(g, \mathbf{x}_{i}g)^{\rho h_R}=(g, \mathbf{x}_{i}^{\rho_{[Y]}}g)^{h_R}=(gh, \mathbf{x}_{i}^{\rho_{[Y]}}gh)$, it follows that $\rho h_R=h_R\rho$.

 Similarly, $ (g, \mathbf{x}_{i}g)^{h_R\tau}=(gh, \mathbf{x}_{i}{gh})^\tau=(kk^{-1}gh, \mathbf{y}_{j}k^{-1}{gh})^\tau
 =(k^{\tau_j}k^{-1}gh, \mathbf{y}_{j}k^{-1}{gh})$ and
 $ (g, \mathbf{x}_{i}g)^{\tau h_R}=(k^{\tau_j}k^{-1}g, \mathbf{y}_{j}k^{-1}g)^{h_R}=(k^{\tau_j}k^{-1}gh, \mathbf{y}_{j}k^{-1}{gh}), $
 so  $\tau h_R=h_R\tau$. Therefore,  $h_R$  commutes with both $\rho$ and $\tau$ which implies $h_R\in \Aut(\mathcal{M} )$.
 \qed

 It is clear that $G_R$ acts regularly on the hypervertex set of Cayley hypermaps of $G$. On the contrary,  a hypermap $\mathcal{M}$ with simple underlying hypergraph is a Cayley hypermap if it has an automorphism group acting regularly on the hypervertex set.

\begin{thm}\label{rightauto}
Let $\mathcal{M}=(D; \rho, \tau)$ be an oriented  hypermap whose underlying hypergraph $\varGamma$ is simple. Then
  $\mathcal{M}$ is a Cayley hypermap if and only if
 $\Aut(\mathcal{M})$ contains a subgroup which induces  regular action on  the hypervertex set of $\G$.
   \end{thm}

\prf Let $\mathcal{G}(\G)$ be the Levi graph of $\G$. Then, $\mathcal{M}$ is the embedding of $\mathcal{G}(\G)$ on an oriented surface.
Assume that $\mathcal{M}$ is a Cayley hypermap of a group $G$, $\mathcal{M}=\CM(G, Y, \rho_{[Y]}, \tau_{[r]})$, where $Y$ is a non-Cayley equivalent subset of $G$ and $|Y|=r$. Lemma~\ref{rightregular} assures the regular action of $G_R$ on $V(\G)$.

On the contrary, assume that $\Aut(\mathcal{M})$ contains a subgroup $G$ which induces regular action, still denoted by $G$, on $V(\G)$.
Then, each hypervertex is adjacent with  the same number of hyperedges in $\mathcal{G}(\G)$, say $d$, and
we can label the hypervertices by the group elements of $G$.

 Choose one hypervertex, say $u$, and label it by $1_G$. Then, a hypervertex $v$ is labelled by $g$ for some $g\in G$ if $v=u^g$. The uniqueness of $g$ is assured by the regular action of $G$ on hypervertices. Now, hyperedges can  be labelled by subsets of $G$. For example, a hyperedge $e$ is labelled by $\{g_1, g_2, \ldots, g_m\}\subseteq G$ if and only if $e$ is adjacent with $u^{g_1}, u^{g_2}, \ldots, u^{g_m}$ in $\mathcal{G}(\G)$. And so we denote $e$ by $\{g_1, g_2, \ldots, g_m\}$, that is, $e := \{g_1, g_2, \ldots, g_m\}$. In this way, $D$ can be written as $$D=\{(g, \mathbf{x})\ |\ g\in \mathbf{x},~ \mathbf{x}\subseteq G~\hbox{and}~\mathbf{x}\in E(\G)\}.$$ Take an arc $(g, \mathbf{x})$ and let $\mathbf{x}=\{g_1=g, g_2, \ldots, g_m\}$. In $\G$, this arc is
$$(u^g, \{u^{g_1}, u^{g_2}, \ldots, u^{g_m}\}).$$
For each $h\in G$,  $$(g, \mathbf{x})^h=(u^g, \{u^{g_1}, u^{g_2}, \ldots, u^{g_m}\})^h=(u^{gh}, \{u^{g_1h}, u^{g_2h}, \ldots, u^{g_mh}\}) =(gh, \mathbf{x}h). $$

Because the cycles of $\rho$ are in one-to-one correspondence with the  hypervertices of $\mathcal{M}$, assume the cycle corresponding with $1_G$ is
$((1_G, \mathbf{x}_1), (1_G, \mathbf{x}_2), \ldots, (1_G, \mathbf{x}_d))$. Then, the cycle of $\rho$  corresponding with $g$ is $((g, \mathbf{x}_1g), (g, \mathbf{x}_2g), \ldots, (g, \mathbf{x}_dg))$ because the actions of $\rho$ and $g$  on $D$ are commutative. That is to say, in $\mathcal{M}$ the cycle $( \mathbf{x}_1g, \mathbf{x}_2g, \ldots,  \mathbf{x}_dg)$ of hyperedges around $g$ is determined by the cycle of hyperedges $( \mathbf{x}_1, \mathbf{x}_2, \ldots,  \mathbf{x}_d)$ around $1_G$. So set $\rho_{[Y]}=( \mathbf{x}_1, \mathbf{x}_2, \ldots,  \mathbf{x}_d)$.

   The cycles of $\tau$ are in one-to-one correspondence with the hyperedges of $\mathcal{M}$. Assume $$((g_1, \mathbf{x}), (g_2, \mathbf{x}), \ldots, (g_m, \mathbf{x}))$$ is a cycle of $\tau$ which corresponds with the hyperedge $\mathbf{x}=\{g_1, g_2, \ldots, g_m\}$ and respects the orientation of $\mathcal{M}$.
       Because the actions of each $h\in G$ and $\tau$ are commutative,
    $$((g_1h, \mathbf{x}h), (g_2h, \mathbf{x}h), \ldots, (g_mh, \mathbf{x}h))$$is a cycle of $\tau$ as well. That is to say, in $\mathcal{M}$ the cycle $(g_1h, g_2h, \ldots, g_mh)$
         of hypervertices around $\mathbf{x}h$ is determined by the cycle of hypervertices $(g_1, g_2, \ldots, g_m)$ around $\mathbf{x}$.

 Consider the action of $G$ on hyperedges and assume this action has $r$ orbits. Then, each orbit contains at least one hyperedge adjacent with $1_G$. In fact, if $\mathbf{x}=\{g_1, g_2, \ldots, g_m\}$ is a hyperedge, then $\mathbf{x}g_i^{-1}=\{g_1g_i^{-1}, g_2g_i^{-1}, \ldots, 1_G, \ldots, g_mg_i^{-1}\}$ is a hyperedge in the orbit of $\mathbf{x}$ for each $1\leq i\leq m$. So, $\mathbf{x}g_i^{-1}$ is adjacent with $1_G$ and
 $(1_G, \mathbf{x}g_i^{-1})$ is an arc.
  Choose one hyperedge adjacent with $1_G$ from each orbit, say $\mathbf{y}_1, \mathbf{y}_2, \ldots, \mathbf{y}_r$, then we have
  $r$ arcs $(1_G, \mathbf{y}_1), (1_G, \mathbf{y}_2), \ldots, (1_G, \mathbf{y}_r)$.
  Set $$Y=\{\mathbf{y}_1, \mathbf{y}_2, \ldots, \mathbf{y}_r\},$$ then $Y$ is a non-Cayley equivalent subset of $G$ and
$[Y]=\{\mathbf{x}_1, \mathbf{x}_2, \ldots, \mathbf{x}_d\}$. Let
  $\widetilde{\G}=\CH(G, Y)$.

   Assume in $\mathcal{M}$, the cycle of hypervertices adjacent with $\mathbf{y}_i$ is $$(1_G, a_{i2},  \ldots, a_{is_i}),$$ then
   define $\tau_i=(1_G, a_{i2}, \ldots, a_{is_i}), 1\leq i\leq r$.
   We claim that $\tau_i$ is an ideal cycle.
   Suppose $\mathbf{y}_i=\mathbf{y}_ih$, then $h\in \mathbf{y}_i$, say $h=a_{ij}$ for some $2\leq j\leq s_i$. As mentioned above, the permutation of $\mathbf{y}_ia_{ij}$ is $(a_{ij}, a_{i2}a_{ij}, \ldots, a_{is_i}a_{ij})$. Look at the arc $(a_{ij}, \mathbf{y}_i)=(a_{ij}, \mathbf{y}_ia_{ij})$. On one hand, $(a_{ij}, \mathbf{y}_i)^\tau=(a_{i,j+1}, \mathbf{y}_i)$. And on the other hand, $(a_{ij}, \mathbf{y}_ia_{ij})^\tau=(a_{i2}a_{ij}, \mathbf{y}_ia_{ij})$. Because $\tau$ is a well-defined permutation on arc set of $\mathcal{M}$, we have $a_{i,j+1}=a_{i2}a_{ij}$. Similarly, one can get $$(1_G, a_{i2}, \ldots, a_{is_i})=(a_{ij}, a_{i2}a_{ij}, \ldots, a_{is_i}a_{ij}),$$ so $\tau_i$ is an ideal cycle.
Note that
 $\rho_{[Y]}=(\mathbf{x}_1, \mathbf{x}_2, \ldots, \mathbf{x}_d)$. Let $\tau_{[r]}=\{\tau_i \di i \in [r]\}$.
 Then
   $\mathcal{M}$ is isomorphic to   $\CM(G, Y, \rho_{[Y]}, \tau_{[r]}) $.
\qed

Let $G$ be a group,  $Y$ be a non-Cayley equivalent hyperset and $\G = \CH(G , Y)$
be the corresponding Cayley hypergraph.
 As stated in Corollary~\ref{normalizer}, $G_R\rtimes \Aut(G, [Y])\le \Aut(\G)$. Take a subset $\mathbf{x}=\{x_1, x_2, \ldots, x_m\}$ of $G$ and a cycle $\t=(g_1, g_2, \ldots, g_k)$ of elements in $ G$. Then $\sigma\in \Aut(G)$ naturally induces an action on $\mathbf{x}$ and on $\t$ as $\mathbf{x}^\sigma=\{x_1^\sigma, x_2^\sigma, \ldots, x_m^\sigma\}$ and
 $\t^\sigma=(g_1^\sigma, g_2^\sigma, \ldots, g_k^\sigma)$, respectively. Moreover, for $a\in G$, set $\t_a=(g_1a, g_2a, \ldots, g_ka)$.
   Let
 $\mathcal{M}=\CM(G, Y, \rho_{[Y]}, \tau_{[r]})$ be a Cayley hypermap with underlying hypergraph $\G$.
 Because $\sigma\in \Aut(G, [Y])$ fixes $[Y]$, it induces a natural action $\tilde{\sigma}$ on arcs of $\G$ as $(g, \mathbf{x}g)^{\tilde{\sigma}}=(g^\sigma, \mathbf{x}^\sigma g^\sigma)$ for any $g\in G, \mathbf{x}\in [Y]$.  It is clear that $\tilde{\sigma}$ is a permutation of $D(\G)$. Hence $\sigma$ can be viewed as a permutation of $D(\G)$.
 As stated in Theorem~\ref{automap},  $\Aut(G, [Y])$ can be viewed as a subgroup of $ \Aut(\mathcal{M})$ under some conditions.

 \begin{thm}\label{automap}
Let $G$ be a group,  $Y=\{\mathbf{y}_1, \mathbf{y}_2, \ldots, \mathbf{y}_r\}$ be a non-Cayley equivalent hyperset of $G$
 and $\mathcal{M}=\CM(G, Y, \rho_{[Y]}, \tau_{[r]})$. Then $G_R\rtimes \Aut(G, [Y])\le \Aut(\mathcal{M})$ if and only if  the following two conditions are satisfied.
\begin{enumerate}
\item [\rm (1)]
For each $\sigma\in \Aut(G, [Y])$ and any $\mathbf{x} \in [Y]$,   $\mathbf{x}^{\sigma\rho_{[Y]}}=\mathbf{x}^{\rho_{[Y]}\sigma}$.

\item [\rm (2)]
For each $\sigma\in \Aut(G, [Y])$ and any  $i\in[r]$, suppose that $\mathbf{y}_i^\sigma=\mathbf{y}_ja$ for some $j \in [r]$ and $a^{-1}\in \mathbf{y}_j$. Then $\t_i^\sigma=\t_{ja}$.
\end{enumerate}
   \end{thm}

\prf
Assume that  $[Y]=\{\mathbf{x}_{1}, \mathbf{x}_{2}, \ldots, \mathbf{x}_{d}\}$,  then
$D(\G)=\{(g, \mathbf{x}_{i}g)\ |\  1\leq i\leq d,  g\in G\}.$

 By Lemma \ref{rightregular}, we have $G_R\le \Aut(\mathcal{M})$. So it suffices to consider $\Aut(G, [Y])$ as a subgroup of $\Aut(\mathcal{M})$. For  $\sigma\in \Aut(G, [Y])$, define $\sigma$ acting on $D(\G)$ by $(g, \mathbf{x}{g})^{\sigma}=(g^\sigma, \mathbf{x}^\sigma g^\sigma)$ for each $\mathbf{x}\in [Y]$ and any $g\in G$. So  $\sigma\in \Aut(\mathcal{M})$ if and only if $\sigma\rho=\rho\sigma$ and $\sigma\tau=\tau\sigma$ with $\rho$ and $\t$ defined in Definition \ref{Cayleyhmap}.

(1) For any $ g\in G$, since
 $$ (g, \mathbf{x}{g})^{\sigma\rho}=(g^\sigma, \mathbf{x}^\sigma g^\sigma)^\rho=
 (g^\sigma, (\mathbf{x}^\sigma)^{\rho_{[Y]}}g^\sigma)$$
 and $$(g, \mathbf{x}{g})^{\rho\sigma}=(g, \mathbf{x}^{\rho_{[Y]}} g)^{\sigma}=
 (g^\sigma, (\mathbf{x}^{\rho_{[Y]}})^\sigma g^\sigma), $$  we have $(g, \mathbf{x}{g})^{\sigma\rho}=(g, \mathbf{x}{g})^{\rho\sigma}$ if and only if $(\mathbf{x}^\sigma)^{\rho_{[Y]}}
 =(\mathbf{x}^{\rho_{[Y]}})^\sigma$. By the arbitrary of $\mathbf{x}$ and $g$, we have  $\sigma\rho=\rho\sigma$ if and only if $(\mathbf{x}^\sigma)^{\rho_{[Y]}}
 =(\mathbf{x}^{\rho_{[Y]}})^\sigma$ for $\mathbf{x} \in [Y]$.

(2) Suppose that $g\in \mathbf{y}_i,h \in G$. Recall that
 $$(gh, \mathbf{y}_ih)^\tau=(g^{\tau_i}h, \mathbf{y}_ih)$$
 in Definition \ref{Cayleyhmap}. Then
 $$(gh, \mathbf{y}_ih)^{\tau\sigma}=(g^{\tau_i}h, \mathbf{y}_ih)^{\sigma}
=(g^{\tau_i\sigma}h^\sigma, \mathbf{y}_jah^\sigma)$$
and $$(gh, \mathbf{y}_ih)^{\sigma\tau}=(g^\sigma h^\sigma, \mathbf{y}_ja h^\sigma)^{\tau}=(u^{\t_j}ah^\sigma, \mathbf{y}_ja h^\sigma)$$ by Definition \ref{Cayleyhmap}, where $g^\sigma=ua \hbox{~and~} u\in \mathbf{y}_j$.
So $(gh, \mathbf{y}_ih)^{\tau\sigma}=(gh, \mathbf{y}_ih)^{\sigma\tau}$ if and only if $g^{\tau_i\sigma}=u^{\t_j}a$ which implies that $\t_i^\sigma=\t_{ja}$ for the arbitrary of $g\in \mathbf{y}_i$.
Therefore, $\sigma\tau=\tau\sigma$ if and only if $\t_{i}^\sigma=\t_{ja}$ for any $i\in [r]$.

Thus we have proofed  $\sigma\in \Aut(\mathcal{M})$ if and only if $(\mathbf{x}^\sigma)^{\rho_{[Y]}}
 =(\mathbf{x}^{\rho_{[Y]}})^\sigma$ for $\mathbf{x} \in [Y]$ and $\t_{i}^\sigma=\t_{ja}$ for any $i\in [r]$.  The theorem follows by the arbitrary of $\sigma$.
\qed

In Example~\ref{examfano},
$\tau_1=(0, 1, 3)$.
In Case 1,
$\rho_{[Y]}=(\{0, 1, 3\},  \{0, 2, 6\},  \{0, 4, 5\})$. Let $\sigma_1\in \Aut(\mathbb{Z}_7)$ which sends $i$ to $2i$ (modular $7$).
A routine check shows that $\sigma_1$ is an automorphism of the corresponding Cayley hypermap which satisfies the conditions in Theorem~\ref{automap}. Therefore, the order of the automorphism group is at least $21$. Because there are $21$ darts, the Cayley hypermap is regular in Case 1. Similarly, take $\sigma_2\in \Aut(\mathbf{Z}_7)$ which sends $i$ to $4i$, then the Cayley hypermap is regular as well in Case 2.

Let $\MM = \CM(G, Y,\rho_{[Y]}, \t_{[r]})$ be a Cayley hypermap, $|Y|=r, [Y]=\{\x_1,\x_2,\ldots,\x_d\}$ and $\rho_{[Y]}=(\x_1,\x_2,\ldots,\x_d)$.
Generally, no further information is known about $\Aut(\mathcal{M})$ except for $G_R\le \Aut(\mathcal{M})$.
 Suppose that there exists a hypermap automorphism $\zeta$ which is not in $G_R$. Because $\zeta$ commutes with $\rho$, it is enough to look at  its action on one $(1_G, \x_i)$ for some $ 1\leq i\leq d$, say $(1_G, \x_d)$.
 Then, $(1_G, \x_d)^\zeta=(g, \x_ig)$ for some $g\in G$ and $i\in [d-1]$. Consequently, there exists a hypermap automorphism which is not in $G_R$ if and only if there exists a hypermap automorphism $\eta=\zeta g^{-1}_R$ such that $(1_G, \x_d)^\eta=(1_G, \x_i)$ for some $i\in [d-1]$.
 Denote the stabilizer of $1_G$ in $\Aut(\MM)$ by
$$\Aut(\MM)_{1_G}=\{\a \in \Aut(\MM) \di (1_G, \x_d)^\a=(1_G, \x_j)~\mbox{for some }~ j \in [d]\},$$ and
let $$\Rot(\MM)=\{j \di (1_G, \x_d)^\a
=(1_G, \x_j),~\a \in \Aut(\MM)_{1_G}\}.$$
Since $\Aut(\MM)$ commutes with $\langle \rho, \tau\rangle$, there is only one element $\a$ in $\Aut(\MM)_{1_G}$ satisfying $(1_G, \x_d)^\a=(1_G, \x_j)$ for each $j\in \Rot(\MM)$.
 So, if $\a \in \Aut(\MM)_{1_G}$ and $(1_G, \x_d)^\a=(1_G, \x_j)$, then we  denote $\a$ by $\a_j$. Since $\a_j\rho=\rho\a_j$,  $(1_G, \x_k)^{\a_j}=(1_G, \x_{k+j})$ for any $k\in [d]$.
The automorphism group of $\MM$ can be characterised using $\Aut(\MM)_{1_G}$ and $\Rot(\MM)$ in the following Theorem \ref{aubch}.

\begin{thm}\label{aubch}
Let $\MM = \CM(G, Y,\rho_{[Y]}, \t_{[r]})$ with $|Y|=r$ and $|[Y]|=d$ be a  Cayley hypermap, and let $m$ be the smallest positive integer of $\Rot(\MM)$. Then the following hold.
\begin{enumerate}
\item [\rm {(1)}]
 $\Aut(\MM)_{1_G}$ is a cyclic group of order $\frac{d}{m}$.

\item [\rm {(2)}] $\Aut(\MM)=\Aut(\MM)_{1_G}G_R$ and $|\Aut(\MM)|=|\Aut(\MM)_{1_G}||G_R|$.
\end{enumerate}
\end{thm}

\prf
Let $[Y]=\{\x_1,\x_2,\ldots,\x_d\}$ and $\rho_{[Y]}=(\x_1,\x_2,\ldots,\x_d)$.

(1) We claim that $\Aut(\MM)_{1_G}=\lg \a_m \rg$. Otherwise, there exists $\a_k \in \Aut(\MM)_{1_G}$ such that
$\a_k \not\in \lg \a_m \rg$. Since $m$ is smallest in $\Rot(\MM)$, $k=qm+r$ for some $q>0$ and $0<r<m$. Then
$(1_G, \x_d)^{\a_k\a_m^{-q}}=(1_G, \x_{k-qm})=(1_G, \x_r)$. Hence $r \in \Rot(\MM)$ which contradicts to the choice of $m$.
Obviously, $|\lg \a_m \rg|=|\a_m|=\frac{d}{m}$.

(2) Since $\Aut(\MM)$ and $G_R$ act transitively on $G$, we have $\Aut(\MM)=\Aut(\MM)_{1_G}G_R$.
Note that $G_R$ acts regularly on $G$. Then $\Aut(\MM)_{1_G}\cap G_R$ only contains the identity automorphism. So,  $|\Aut(\MM)|=|\Aut(\MM)_{1_G}||G_R|$.
\qed

\begin{cor}
Let $\MM = \CM(G, Y,\rho_{[Y]}, \t_{[r]})$ with $|Y|=r$ and $|[Y]|=d$ be  a Cayley hypermap. Then  the following hold.
\begin{enumerate}
\item [\rm {(1)}]
$\MM$ is  regular
if and only if $\Aut(\MM)_{1_G}$ is a cyclic group of order $d$.

\item [\rm {(2)}] If $d$ is a prime number, then $\MM$ is regular or $\Aut(\MM)$ is isomorphic to $ G_R$.
\end{enumerate}
\end{cor}

\vskip 0.5cm

Let $G$ be a finite group, $\sigma\in \Aut(G)$ and $|\sigma|=k$. Take $a\in G$ and let $$\mathcal{O}_a =
\{a^\mu \di \mu\in \lg \sigma \rg\}.$$
Assume that $|\mathcal{O}_a|=k$. Let
 $X=\{\lg b\rg \di b\in \mathcal{O}_a\}$.  Then $X$ is inverse closed.
 The equality $\lg a^{\sigma^i}\rg=\lg a^{\sigma^j}\rg$ is equivalent to $\lg a\rg=\lg a^{\sigma^{j-i}}\rg$ for $i,j\in [k]$. So, to make sure that $\G=\CH(G, X)$ is a  $k$-uniform simple Cayley hypergraph, we assume that $\lg a\rg\ne \lg a^{\sigma^{i}}\rg$ for any $i\in [k-1]$. To be convenient, denote $\mathbf{y}_i=\lg a^{\sigma^i}\rg, 1\leq i\leq k$.
 Let $$\rho_X=(\mathbf{y}_1, \mathbf{y}_2, \ldots, \mathbf{y}_{k}),$$
and
$$\tau_i=(1_G, a^{\sigma^i}, (a^2)^{\sigma^i},\ldots, (a^{t-1})^{\sigma^i})\hbox{~for~all~}1\leq i\leq k,$$ where $t$ is the order of $a$. It is easy to check that $\tau_i$ is an ideal permutation of $\mathbf{y}_i$ for each $1\leq i\leq k$.
Then $\MM = \CM(G, X,\rho_{X}, \t_{[k]})$ is a Cayley hypermap.
\begin{thm}\label{th3.14}
Using the notation above, suppose that $\lg \mathcal{O}_a \rg=G$.
Then the hypermap $\MM$ is a regular Cayley hypermap with the simple underlying hypergraph $\G$ and $\Aut(\MM)=\lg \sigma\rg G_R$.
\end{thm}
\prf
The arc set of $\MM$ is $D=\{(g, \mathbf{y}_ig)\ |\ g\in G, 1\leq i\leq k\}$. Obviously, the underlying hypergraph of $\MM$ is simple. Let $\mathbf{y}_{k+1} =\mathbf{y}_1$ and $\mathbf{y}_{k+2} =\mathbf{y}_2$.

Note that $\sigma$ is a permutation on $G$,  and so the $\sigma$-action on $G$ can be extended to $D$ by defining
$$(g, \mathbf{y}_ig)^{\sigma}=(g^\sigma, \mathbf{y}_i^\sigma g^\sigma) \hbox{~for~all~} i\in [k] \hbox{~and~} g\in G.$$
Thus $(g, \mathbf{y}_ig)^{\sigma}=
(g^\sigma, \mathbf{y}_{i+1}g^\sigma) \hbox{~for~} i\in [k] \hbox{~and~} g\in G$. Now let $i\in [k]$ and $g\in G$. Then we get
$(g, \mathbf{y}_ig)^{\tilde{\sigma}\rho}=(g^\sigma, \mathbf{y}_{i+2}g^\sigma)$ and $$(g, \mathbf{y}_ig)^{\rho\sigma}=(g, \mathbf{y}_{i+1}g)^{\sigma}=(g^\sigma, \mathbf{y}_{i+2}g^\sigma) .$$ Thus $\sigma\rho=\rho\sigma$.
Similarly, we have $(g, \mathbf{y}_ig)^{\sigma\tau}=(a^{\sigma^{i+1}}g^\sigma, \mathbf{y}_{i+1}g^\sigma)$ and $$(g, \mathbf{y}_ig)^{\tau\sigma}=(a^{\sigma^i}g, \mathbf{y}_{i}g)^{\sigma}=(a^{\sigma^{i+1}}g^\sigma, \mathbf{y}_{i+1}g^\sigma).$$ So $\sigma\tau=\tau\sigma$. These equalities imply that $\sigma\in \Aut(\MM)$. Especially, $(1_G, \mathbf{y}_i)^{\sigma}=(1_G, \mathbf{y}_{i+1})$.
According to Theorem~\ref{aubch}, $\Aut(\MM)_{1_G}=\langle \sigma\rangle$. So $\Aut(\MM)=\lg \sigma\rg G_R$ and $\MM$ is a regular Cayley hypermap.
\qed
\begin{exam} Let $G=\overbrace{\mathbb{Z}_p\times \mathbb{Z}_p\times \cdots \times \mathbb{Z}_p}^k$ be an elementary abelian group, where $p$ is an odd prime number. Assume that the $i$-th $\mathbb{Z}_p$ is generated by $a_i$ for each $1\leq i\leq k$. Define an automorphism  $\sigma$ of $G$ as $a_i^\sigma=a_{i+1}$ for $i\in [k-1]$ and $a_k^\sigma=a_{1}$.
It is easy to check that $|\sigma|=k$.
Let $\mathbf{y}_i=\langle a_i\rangle,  X=\{\mathbf{y}_i \di 1\leq i\leq k\}, \rho_X=(\mathbf{y}_1, \mathbf{y}_2, \ldots, \mathbf{y}_k)$ and $$\tau_i=(1_G, a_i, a_i^2, \ldots, a_i^{p-1})~\hbox{for~each~} i\in [k].$$ Then, according to Theorem~{\rm \ref{th3.14}}, $\CM(G, X, \rho_X, \tau_{[k]})$ is a regular Cayley hypermap.
\end{exam}

\section{Further research}

In this last section we propose some questions for further
research.

1. For which groups $G$ is there a regular Cayley hypermap (with 3-uniform underlying hypergraphs)? Classify regular Cayley hypermaps for cyclic groups, dihedral groups, elementary abelian $p$-groups and simple groups.

2. Which regular hypermaps of small genus or prime-power vertices are Cayley hypermaps?

3. Characterize normal Cayley hypermaps whose full automorphism group contains a normal regular subgroup acting on vertices.

4. What is the theory of Cayley hypermaps with non-orientable surfaces?

\section*{Acknowledgements}
The authors thank the referees for their helpful comments and suggestions. This work  was
supported by  the National Natural Science Foundation of China (No. 12101535) and NSFS (No. ZR2020MA044).

\section*{Data Availability}

Data sharing is not applicable to this article as no datasets were generated or analyzed
during the current study.

\section*{Declarations}

{\bf Conflict of interest} The authors declare that they have no conflict of interest.

\end{document}